\documentclass[12pt]{article}

\usepackage{amsmath}    
\usepackage{amssymb}
\usepackage{graphicx}   
\usepackage{verbatim}   
\usepackage{color}      

\usepackage{hyperref}   

\usepackage{enumerate}
\usepackage{mathrsfs}
\usepackage{empheq}
\usepackage{bbold}

\usepackage[final]{showlabels}

\usepackage{amsthm}

\newtheorem{lemma}{Lemma}[section]
\newtheorem{corollary}{Corollary}[section]
\newtheorem{theorem}{Theorem}[section]

\newcommand{\del}{\partial}
\renewcommand{\theta}{\vartheta}
\renewcommand{\phi}{\varphi}
\newcommand{\vecc}[2]{\left ( \begin{array}{c}#1\\#2\\ \end{array}\right )}
\newcommand{\veccc}[3]{\left ( \begin{array}{c}#1\\#2\\#3\\ \end{array}\right )}

\newcommand{\dd}{\mathrm{d}}

\newcommand{\id}{\mathbb{1}}

\renewcommand{\vec}{\mathbf}

\newcommand{\vecalpha}{\boldsymbol\alpha}

\newcommand{\vecxi}{\boldsymbol\xi}

\newcommand{\const}{\mathrm{const}}

\newcommand{\vecdelta}{\boldsymbol\delta}

\renewcommand{\title}{The Active Flux scheme for nonlinear problems}

\newcommand{\authorOne}{Wasilij Barsukow\footnote{Institute for Mathematics, Zurich University, 8057 Zurich, Switzerland}}

\begin{document}

\begin{center} \Large
\title

\vspace{1cm}

\date{}
\normalsize

\authorOne
\end{center}

\begin{abstract}

The Active Flux scheme is a finite volume scheme with additional point values distributed along the cell boundary. It is third order accurate and does not require a Riemann solver. Instead, given a reconstruction, the initial value problem at the location of the point value is solved. The intercell flux is then obtained from the evolved values along the cell boundary by quadrature. Whereas for linear problems an exact evolution operator is available, for nonlinear problems one needs to resort to approximate evolution operators. This paper presents such approximate operators for nonlinear hyperbolic systems in one dimension and nonlinear scalar equations in multiple spatial dimensions. They are obtained by estimating the wave speeds to sufficient order of accuracy. Additionally, an entropy fix is introduced and a new limiting strategy is proposed. The abilities of the scheme are assessed on a variety of smooth and discontinuous setups.



\end{abstract}

\section{Introduction}

Hyperbolic $m\times m$ systems of conservation laws in $d$ spatial dimensions have the form
\begin{align}
 \del_t q + \nabla \cdot \vec f(q) &= 0 & q: \mathbb R^+_0 \times \mathbb R^d \to \mathbb R^m \label{eq:systemconslaw}
\end{align}
The function $\vec f$ is called the flux. Exact solutions of these equations are unavailable in general, and one needs to resort to numerical methods. 

Cell based methods consider the computational domain to be partitioned into cells. A certain number of discrete degrees of freedom are associated with every cell: e.g. finite volume methods store the average of the dependent variable and spectral/Galerkin methods store coefficients of a decomposition in some basis. 

In order to evolve the cell average, \emph{finite volume methods} require the know\-ledge of the flux through the intercell boundary (see section \ref{ssec:finvol} for a derivation). This flux cannot generally be approximated by a symmetric average of fluxes associated to the values in the two adjacent cells, because this results in an unstable method. Instead, the fact that hyperbolic PDEs have certain preferred directions of information propagation needs to be reflected in the numerical method. The choice of the numerical flux as an asymmetric average of the neighbouring values is referred to as \emph{upwinding}. It has been suggested in \cite{godunov59difference} to use an exact short-time solution as a building block in order to find a numerical flux that leads to a stable scheme. First, a piecewise polynomial function is found, such that it is continuous in every cell and its average agrees with the given average. The discontinuities at cell interfaces present so-called Riemann Problems, which then are solved over a time interval that does not allow them to interact. The exact flux at the location of the cell interface then is used as a numerical flux in the finite volume method. To save computation time, an approximate solution of the Riemann Problem can be used, see e.g. \cite{roe81,harten83,jin95} as well as \cite{leveque02,toro09} for more details. Higher order of accuracy is achieved by widening the stencil (\cite{vanleer77,colella84,titarev02}).

Galerkin methods represent the numerical solution in e.g. a polynomial basis. Every basis coefficient is then evolved using the weak formulation of \eqref{eq:systemconslaw}. Again, for hyperbolic equations this requires modification in order to achieve a stable method, one of which is the \emph{discontinuous Galerkin method} (\cite{cockburn98}, but see also \cite{brooks82}). The basis functions are piecewise polynomial, and in order to deal with the jumps across cell interfaces a Riemann solver is invoked. Higher order of accuracy is achieved by retaining more coefficients of the basis decomposition. 

Even in one spatial dimension, conservation laws \eqref{eq:systemconslaw} therefore pose a number of challenges to numerical methods. This does not only include the necessity of upwinding. It is also known that continuous solutions do not generally exist for all times, and thus numerical methods need to be designed in such a way that they can capture discontinuities (weak solutions). Weak solutions in one spatial dimension only become unique upon additional conditions (entropy conditions), and numerical methods need to fulfill a discrete counterpart of these conditions (entropy stability, see e.g. \cite{tadmor2003entropy}). These aspects have been subject of numerous investigations, see e.g. \cite{leveque02} for an introduction.

In multiple spatial dimensions, systems of conservation laws have a rich phenomenology which is absent in the one-dimensional case. In the context of the Euler equations these are vortices (e.g. created by Kelvin-Helmholtz instabilities), multi-dimensional shock interactions, the low Mach number/in\-compressible limit and many more. The easiest way of extending a one-dimen\-sional numerical method to multiple dimensions is directional splitting, i.e. the problem is replaced by a number of one-dimen\-sional problems. This, however, has been demonstrated to require excessive grid refinement in order to capture truly multi-dimensional features even for systems much simpler than the Euler equations (e.g. \cite{morton01,guillard04,barsukow17a,barsukow17}). It has been found that numerical methods should reflect essential properties of the solution at discrete level in order to avoid expensive grid refinement. Such methods are called \emph{structure preserving}. So far, modifications of existing schemes have been suggested, but it is largely unexplored how such schemes can be derived from first principles.

The \emph{Active Flux} scheme is a new scheme (\cite{eymann13}, an extension of \cite{vanleer77}) that combines a finite volume scheme with additional, independently evolved degrees of freedom which are interpreted as point values. These point values are located at cell boundaries and Active Flux thus uses a continuous reconstruction. This is a major difference to finite volume schemes. The Active Flux scheme is nevertheless able to resolve shocks (which are approximated by steep gradients), as can be seen below. The reconstruction is parabolic and therefore Active Flux is third order accurate. It has been shown in \cite{barsukow18activeflux} for the equations of linear acoustics that the scheme is vorticity preserving without any fix which makes it a good candidate as a structure preserving method for more complicated multi-dimensional problems.

So far, the Active Flux scheme has been studied in great detail for linear equations (\cite{vanleer77,eymann13,barsukow18activeflux}). As explained in section \ref{sec:af}, the essential ingredient is an approximate solution operator for the initial value problem which is used to update the pointwise degrees of freedom. For linear equations, the point values can be updated using an exact evolution operator. 

For nonlinear problems, an approximate evolution operator is required. By exploiting special properties the Active Flux scheme has been applied to Burgers' equation (\cite{eymann11,eymann11a,roe17a}) and Euler equations (\cite{eymann11,fan17,maeng17,kerkmann18}). Some of these extensions lose the order of convergence when applied to nonlinear equations. The approximate evolution operator has to be of sufficiently high order, e.g. local linearization as used in \cite{eymann11} is not sufficient to yield an overall third order scheme. In \cite{kerkmann18}, a solution operator based on the Cauchy-Kovalevskaya/\-Lax-Wen\-droff procedure has been suggested which can be applied to general nonlinear hyperbolic systems in one spatial dimension and has shown the correct order of convergence when applied to one-dimensional Euler equations in practice. However, both the procedure of \cite{kerkmann18} itself and the evaluation of the higher order spatial derivatives can be rather complicated. In particular, the derivatives are required at locations where the reconstruction is not differentiable.

The aim of this paper is to provide a simpler solution operator that allows to apply the Active Flux scheme to a large class of hyperbolic conservation laws. The general idea is to keep the structure of a characteristic-based (or in multi-d characteristic-cone-based, see \cite{barsukow18activeflux}) evolution operator but to estimate carefully the wave speeds -- which are not constant in the nonlinear case. This also includes an estimate on whether a shock has occurred by self-steepening. In section \ref{sec:scalar} such approximate evolution operators are provided for scalar conservation laws in one and several spatial dimensions, and in section \ref{sec:systems} -- for hyperbolic systems of conservation laws in one spatial dimension. 
This leads to algorithms significantly different from the scalar case. This paper is the first part of a sequence of papers devoted to the application of the Active Flux scheme to nonlinear problems. The case of multi-dimensional systems shall form the content of a forthcoming work. Although the examples presented here are all computed on Cartesian grids, the approximate evolution operators can be immediately applied to unstructured grids. A detailed experimental study concerning unstructured grids, however, is subject of future work.

Section \ref{sec:limiting} is describing a limiting procedure. As the Active Flux scheme is of higher order, spurious oscillations can appear. The continuous reconstructions employed in the Active Flux scheme do not allow to make immediate use of the same limiting strategies as in the case of usual finite volume schemes. Several limiters for Active Flux have been suggested in the literature: in \cite{roe15,kerkmann18} the parabolic reconstruction inside a cell is replaced by several parabolae joined in a continuous and monotone manner. This, however, might add implementational and computational complexity. The same is true for the hyperbolic reconstruction considered in \cite{kerkmann18} as its parameters cannot be computed analytically. Additionally, discontinuous reconstructions have been considered in \cite{eymann11,eymann13a,kerkmann18} as limiting strategies. However, favorable properties have been deduced from the continuous reconstruction in \cite{barsukow18activeflux,barsukow17a}, and these discontinuous limiting strategies violate the principle of Active Flux and move it again closer to usual finite volume schemes. This shows the need for a simple limiter that keeps the reconstruction continuous. Here, such a limiter is presented -- it is optimally monotone (see Theorem \ref{thm:powerlawlimiting} for precise statement) and is at the same time computationally efficient.

The paper thus is organized as follows: Approximate evolution operators are presented in section \ref{sec:scalar} for scalar conservation laws and in section \ref{sec:systems} for systems. Limiting is discussed in section \ref{sec:limiting} and numerical examples for problems in one and two spatial dimensions are shown in section \ref{sec:numerical}.

\section{The Active Flux scheme} \label{sec:af}

\subsection{Finite volume scheme}\label{ssec:finvol}
Consider the computational domain to be divided into (polyhedral) computational cells $\mathcal C \subset \mathbb R^d$ and discretize the time into points $t^n$, $n\in \mathbb N_0$ separated by (not necessarily equal) time steps $\Delta t$. Recall that in order to solve \eqref{eq:systemconslaw}, finite volume schemes use cell averages\footnote{Boldface symbols are reserved for elements of vector spaces of the same dimension $d$ as the space, if $d\geq 2$. Indices never denote derivatives.}
\begin{align}
 \bar q_{\mathcal C} = \frac{1}{|\mathcal C|} \int_{\mathcal C} \dd \vec x \, q(t, \vec x)
\end{align}
as discrete degrees of freedom. The cell average is updated in time using fluxes $f_e$ through cell boundaries (edges $e$):
\begin{align}
 \frac{\bar q_{\mathcal C}^{n+1} - \bar q_{\mathcal C}^n}{\Delta t} + \sum_{e \subset \del \mathcal C} \frac{|e|}{|\mathcal C|} f_e &= 0 \label{eq:finitevolume}
\end{align}

Here, $\bar q_{\mathcal C}^n$ denotes the value of the average at time $t^n$.

Applying Gauss' law to \eqref{eq:systemconslaw}, $f_e$ can be given the interpretation of approximating
\begin{align}
 f_e \simeq \frac{1}{\Delta t} \int_{t^n}^{t^{n+1}} \dd t\, \frac{1}{|e|}\int_e \dd x \vec n_e \cdot \vec f(q) \label{eq:numflux}
\end{align}
with $\vec n_e$ the outward normal to edge $e$.

In the one-dimensional ($d=1$) case the cells are indexed by a finite subset of the integers. Then $\bar q_i$ denotes the averages in cell $\mathcal C_i = [x_{i-\frac12}, x_{i+\frac12}]$, which is centered around $x_i$. The size $|C_i|$ of a cell is for simplicity denoted by $\Delta x$. Most of the results remain valid when the size of the cells varies smoothly.

\subsection{Pointwise degrees of freedom}

\begin{figure}
 \centering
 \includegraphics[width=0.6\textwidth]{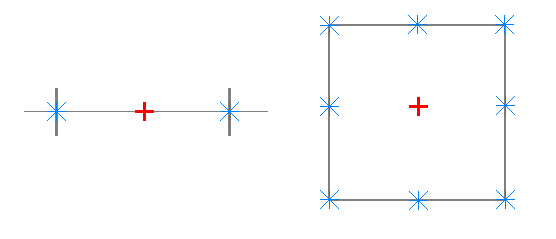}
  \caption{The degrees of freedom used for Active Flux. Stars indicate the location of point values, and the cross (placed in the center symbolically) refers to the cell average. \emph{Left}: One spatial dimensions. \emph{Right}: Two spatial dimensions.}
 \label{fig:dofactiveflux}
\end{figure}

The Active Flux scheme is an extension of the finite volume scheme \eqref{eq:finitevolume} for \eqref{eq:systemconslaw}. Active flux uses point values $q_\vec x$ located at points $\vec x \in \mathbb R^d$ along the cell boundary as additional discrete degrees of freedom (recall that indices never denote derivatives in this paper). They approximate the value $q(t, \vec x)$. So far, the following choices have been considered in the literature (see also Figure \ref{fig:dofactiveflux}):
\begin{itemize}
\item In one spatial dimension the point values are located at cell boundaries $x_{i+\frac12}$, $i\in\mathbb Z$ and are thus rather denoted by $q_{i+\frac12}$. 
\item In two spatial dimensions so far (in \cite{eymann13,barsukow18activeflux}) the locations of the pointwise degrees of freedom are chosen to be the endpoints and the midpoints of edges.
\end{itemize}

Active flux is not a staggered-grid finite volume method. Staggered grids consider offset grids of averages, each for a different variable; Active Flux stores additional point values of all the variables inside every cell. This approach is closer to Lagrange-basis spectral/Galerkin methods, with the difference that the cell average is retained as one of the degrees of freedom. A particularity of the Active Flux method is the exclusive distribution of the point values along the cell boundary.

For the evolution of the point values at cell boundaries the Active Flux scheme considers an initial value problem: a reconstruction $q_\text{recon}(\vec x)$ plays the role of the initial data, and the evolution in time can be either exact or approximate. This is explained in more detail in the next sections.

Once the time evolution of the point values is known, the numerical flux is obtained using a quadrature of \eqref{eq:numflux} in time and along the edge. In order to obtain a third order scheme it is necessary to compute the point values also at half the time step (see \cite{barsukow18activeflux} for further implementation details). Only the update of the average needs to be conservative, there is no notion of a conservative update for a point value.

\subsection{Reconstruction} \label{ssec:introreconstruction}

The reconstruction is interpolating the point values and the average:
\begin{align} 
 q_\text{recon}(\vec x) &= q_\vec x \quad \forall\text{ locations $\vec x$ of the pointwise degrees of freedom} \label{eq:interpolpropertyrecon}\\
 \frac{1}{|\mathcal C|}\int_\mathcal C \dd \vec x \, q_\text{recon}(\vec x) &= \bar q_{\mathcal C} \label{eq:conservpropertyrecon}
\end{align}
The reconstruction thus is conservative. The difference to reconstructions in the context of finite volume schemes is the fact that the reconstruction is continuous at the locations of the pointwise degrees of freedom. Additionally, the choices used for the reconstruction so far in the literature were such that the reconstruction is continuous \emph{everywhere} in the computational domain. These particular choices are briefly reviewed next:
\begin{itemize}
\item In one spatial dimension, the reconstruction is chosen piecewise parabolic in \cite{vanleer77}. This is a natural choice, as \eqref{eq:interpolpropertyrecon}--\eqref{eq:conservpropertyrecon} amount to three conditions in each cell. It reads
\begin{align}
 q_\text{recon}(x) &= -3 (2 \bar q_i - q_{i-\frac12} - q_{i+\frac12}) \frac{(x-x_i)^2}{\Delta x^2} \\&+ (q_{i+\frac12} - q_{i-\frac12}) \frac{x-x_i}{\Delta x} + \frac{6 \bar q_i - q_{i-\frac12} - q_{i+\frac12}}{4}  \qquad x \in [x_{i-\frac12},x_{i+\frac12}]\label{eq:parabolicrecon}
\end{align}
The reconstruction is continuous everywhere.
\item In two spatial dimensions, in \cite{eymann13,barsukow18activeflux} the pointwise degrees of freedom are placed at endpoints and at the midpoints of every edge. The reconstruction is chosen always to reduce to a parabola along any edge and, as a parabola is uniquely defined by three points, the reconstruction is thus continuous across any edge, and thus everywhere.
\end{itemize}

\subsection{Evolution of pointwise degrees of freedom} \label{ssec:intropointupdate}

Active flux is a time-explicit method and thus subject to a CFL condition
\begin{align}
 \Delta t < \frac{L_\text{min}}{\lambda_\text{max}}
\end{align}
In the following, the time step is chosen based on the maximum value $\lambda_\text{max}$ of the characteristic speed at the location of pointwise degrees of freedom. The shortest length $L_\text{min}$ in the one-dimensional case is the size of the cell, and in the two-dimensional case half the edge length (as there is a pointwise degree of freedom located at its midpoint).

The update procedure of the Active Flux scheme for the point value is the (exact or approximate) solution of the initial value problem at its location. The initial data are given by the reconstruction. When the Active Flux scheme is applied to linear equations (as in \cite{vanleer77,eymann13,barsukow18activeflux}) an exact evolution is easily available. For nonlinear equations it is necessary to devise approximate evolution operators. This is the topic of section \ref{sec:scalar} (for scalar nonlinear conservation laws) and section \ref{sec:systems} (for systems of conservation laws). Here, only a general statement shall be given that concerns the necessary accuracy of an approximate evolution operator. First, the following result is needed:

\begin{lemma} \label{lem:errorcancellation}
 For $f: \mathbb R \times \mathbb R \to \mathbb R$ and $g : \mathbb R \to \mathbb R$, both analytic and $n_1, n_2 \in \mathbb N$, assume
 \begin{align}
  f(x, \Delta x) &= g(x) + \Delta x^{n_1} g^{(n_2)}(x)  + \mathcal O((\Delta x)^{n_1+1}) \qquad \forall x
 \end{align}
 Then
 \begin{align}
  f(x + \Delta x , \Delta x) - f(x, \Delta x) =  g(x + \Delta x) - g(x) + \mathcal O((\Delta x)^{n_1+1})
 \end{align}
\end{lemma}
\begin{proof}
 Expand
 \begin{align}
  f(x + \Delta x , \Delta x) - f(x, \Delta x) &= g(x + \Delta x)  - g(x)\\\nonumber+ \Delta x^{n_1} g^{(n_2)}(x + \Delta x) &- \Delta x^{n_1} g^{(n_2)}(x)  + \mathcal O((\Delta x)^{n_1+1})\\
  &= g(x + \Delta x)  - g(x)\\\nonumber+ \Delta x^{n_1} \Big ( g^{(n_2)}(x) + \mathcal O(\Delta x) \Big ) &- \Delta x^{n_1} g^{(n_2)}(x)  + \mathcal O((\Delta x)^{n_1+1})
 \end{align}
\end{proof}\hfill$\Box$

Recall that $f \in \Theta(g)$ means that asymptotically $ c_1 |g| \leq |f| \leq c_2 |g|$ for some $c_1, c_2 > 0$.

\begin{theorem}\label{thm:order}
 Assume a hyperbolic CFL condition $\Delta x \in \Theta(\Delta t)$ as $\Delta t \to 0$. If the approximate evolution $\tilde q(t, x)$ for fixed $x \in \mathbb R$ approximates the exact solution $q(t, x)$ at least as 
\begin{align}
\tilde q(t, x) = q(t, x) + \mathcal O(t^3)
\end{align}
and the quadrature rules used to approximate \eqref{eq:numflux} yield the exact value up to an error of $\mathcal O(\Delta t^\alpha \Delta x^\beta)$, $\alpha + \beta \geq 3$ then Active Flux formally achieves third order accuracy.
\end{theorem}

\begin{proof}
 Denote by $T_t[q_0]$ the exact evolution operator applied to initial data $q_0$ and evolving them to a time $t$, and by $\tilde T_t[q_0]$ its corresponding approximation.
 
 Assume point values of $q(t^n, x)$ to be used in the reconstruction. Then, because the reconstruction is an interpolation, and taking $x$ to be the location of one of the point values (where the interpolation is exact)
 \begin{align}
 q_\text{recon}^n(x + \delta x) = q(t^n, x) + \mathcal O((\delta x)^\alpha \Delta x^\beta) \qquad \text{with } \alpha + \beta \geq 3 ,\beta \geq 1
 \end{align}
 This statement can also be understood as follows: the interpolation matches the Taylor series $$q(t^n, x) + \del_x q(t^n, x)\delta x + \frac12 \del_x^2 q(t^n, x)(\delta x)^2 + \mathcal O((\delta x)^3)$$ of $q(t^n, x + \delta x)$ in $\delta x$ to sufficiently high powers of $\delta x$. At the same time, the derivatives if $q$ that appear as coefficients in this Taylor series are approximated by finite differences, which carry error terms $\mathcal O(\Delta x^\beta)$, $\beta > 1$. Because they all use the same point values in the approximation, lower order derivatives are approximated better.

 The approximate evolution operator uses initial data from the neighbouring cells at a distance $\mathcal O(\Delta t)$ from some fixed $x$. Therefore
 \begin{align}
 \tilde T_{\Delta t}[q_\text{recon}^n](x) &= T_{\Delta t}[q_\text{recon}^n](x) + \mathcal O(\Delta t^3) 
 \\&=  T_{\Delta t}[q(t^n, \cdot)](x) + \mathcal O(\Delta t^\alpha \Delta x^\beta)   \qquad \text{with } \alpha + \beta \geq 3
 \end{align}

 For the average update, the numerical flux is obtained using a quadrature of \eqref{eq:numflux}, such that the numerical flux differs from the exact one by the same error. Thus, using the assumption of a hyperbolic CFL constraint, lemma \ref{lem:errorcancellation} implies that the leading errors cancel when the fluxes at $x + \Delta x$ and $x$ are subtracted. One is left with
 \begin{align}
  \bar q_i^{n+1} &= \bar q_i^n - \frac{\Delta t}{\Delta x} (f_{i+\frac12}-f_{i-\frac12}) 
  = \bar q_i^n + \text{exact flux difference} + \mathcal O(\Delta t^4)
 \end{align}
 This on total gives a numerical method of third order.
\end{proof}\hfill$\Box$

\subsection{Overview of the algorithm}

The overall algorithm of Active Flux is as follows:

\begin{enumerate}
 \item Given cell averages and point values, compute a reconstruction according to section \ref{ssec:introreconstruction}.
 \item Use the reconstruction as initial data in the update of the point values (section \ref{ssec:intropointupdate}). Approximate evolution operators for scalar nonlinear problems are discussed in section \ref{sec:scalar} and for nonlinear systems in one spatial dimensions in section \ref{sec:systems} below.
 \item Given the updated point values along the cell interfaces, compute the intercell fluxes via quadrature of \eqref{eq:numflux}. Here, a space-time Simpson rule is used.
 \item Update the cell averages via \eqref{eq:finitevolume}.
\end{enumerate}

\section{Scalar nonlinear equations} \label{sec:scalar}

Consider the initial value problem for the following scalar (i.e. $m=1$) conservation law
\begin{align}
 \del_t q + \nabla \cdot \vec f(q) &= 0 & q : \mathbb R^+_0 \times \mathbb R^d &\to \mathbb R \label{eq:scalar}\\ 
 q(0, \vec x) &= q_0(\vec x) & \vec f : \mathbb R  &\to \mathbb R^d
\end{align}

Assume the flux function to be smooth and convex. 

\subsection{Fix-point iteration}

In the absence of shocks \eqref{eq:scalar} can be rewritten as
\begin{align}
 \del_t q + \vec a(q) \cdot \nabla q &= 0
\end{align}
with $\vec a(q) = \del_q \vec f(q)$. The characteristics $\vecxi : \mathbb R^+_0 \to \mathbb R^d$ are straight lines on which the solution is constant. They fulfill
\begin{align}
 \vecxi'(\cdot) &=  \vec a\Big(q(\cdot, \vecxi(\cdot))\Big) &
 \vecxi(t) &= x
\end{align}
The exact solution is found by evaluating the initial data at the footpoint $\hat \vecxi = \vecxi(0)$ of the characteristic
\begin{align}
 q(t, \vec x) = q_0(\hat\vecxi)
\end{align}
as $q$ remains constant along it. This also allows to write
\begin{align}
 \vec x = \hat \vecxi + \vec a(q_0(\hat\vecxi)) t \label{eq:characteristic}
\end{align}

This equation can be solved for $\hat\vecxi$ efficiently using a fixpoint iteration:
\begin{theorem} \label{thm:fixpointiteration}
 $\hat \vecxi^{(n)}$, given recursively by
 \begin{align}
  \hat\vecxi^{(0)} &= \vec x \label{eq:scalariterationsimple1}\\ 
  \hat\vecxi^{(n)} &= \vec x - \vec a(q_0(\hat \vecxi^{(n-1)})) t \qquad n = 1, 2, \ldots \label{eq:scalariterationsimple2}
 \end{align}
 for $t \geq 0$, formally approximates $\hat\vecxi$ to $n$-th order, i.e. $\hat\vecxi^{(n)} = \hat\vecxi + \mathcal O(t^{n+1})$.
\end{theorem}
\begin{proof}
 Define the error $\epsilon^{(n)} \vec d^{(n)} := \hat\vecxi^{(n)} - \hat \vecxi$ with $\|\vec d\| = 1$, $\epsilon^{(n)} \geq 0$ and $\vec A := \vec a \circ q_0$. Then
 \begin{align}
  \hat\vecxi + \epsilon^{(n)} \vec d^{(n)} =\hat\vecxi^{(n)} &\overset{\eqref{eq:scalariterationsimple2}}{=} \vec x - \vec A(\hat \vecxi + \epsilon^{(n-1)} \vec d^{(n-1)})t\\
   &= \vec x - \vec A(\hat \vecxi) t-  \sum_{i = 1}^\infty \vecalpha_i \cdot (\epsilon^{(n-1)})^i \cdot t\\
  \intertext{where $\vecalpha_i = \frac{1}{i!} \nabla_{\vecxi} \vec A \big |_{\hat \vecxi} \cdot \vec d^{(n-1)}$ }
  \epsilon^{(n)} = \|\epsilon^{(n)} \vec d^{(n)}\| &= \Big\|\sum_{i = 1}^\infty \vecalpha_i \cdot (\epsilon^{(n-1)})^i \cdot t\Big\|
 \end{align}
 Obviously $\epsilon^{(0)} \in \mathcal O(t)$. Then by induction, if $\epsilon^{(n-1)} \in \mathcal O(t^{n})$, then for $n \geq 1$ and some constant $C \geq 0$
 \begin{align}
  \epsilon^{(n)} &\leq C \cdot \sum_{i = 1}^\infty (\epsilon^{(n-1)})^i \cdot t\\
  \epsilon^{(n)} &\in \mathcal O( t^{n+1})
 \end{align}
 which proves the assertion.
 \end{proof}\hfill$\Box$
 
This iteration seems related to the Picard iteration, but it is exact for linear problems for \emph{any} initial data after \emph{one} step. Therefore the above iteration is even more powerful than a standard Picard iteration.

In view of Theorem \ref{thm:order}, an evolution operator for the discrete degree of freedom $q_\vec x$ located at $\vec x$ is
\begin{align}
 q_\vec x^{n+1} = q_\text{recon}(\hat \vecxi^{(2)}) \label{eq:evoscalarsimple}
\end{align}
and instead of $q_0$ the reconstruction $q_\text{recon}$ based on values at time $t^n$ would be used in the fixpoint iteration.

\subsection{Comparison to previous results}

Before turning to questions regarding the possible presence of shocks in the solution, compare this evolution to similar approaches available in the literature. Note that \eqref{eq:evoscalarsimple} estimates the speed of the characteristic as
\begin{align}
 \vec a(q_0(\hat \vecxi^{(1)})) = \vec a(q_0(\vec x - \vec a(q_0(\vec x)) t))
\end{align}
Local linearization would correspond to taking the evolution operator $q_\text{recon}(\hat \vecxi^{(1)})$, and thus estimate the characteristic speed simply by $\vec a(q_0(\vec x))$. For the special case of Burgers' equation, in \cite{eymann11} it is suggested to estimate the characteristic speed in one spatial dimension by
\begin{align}
 \frac{1}{2} (q_{i+\frac12} + q_{i-\frac12})
\end{align}
However, this approach does not lead to an increase in the order of convergence (as can be shown by direct computation) and thus is not fundamentally superior to local linearization.

In \cite{roe17a} the exact speed of the characteristic for \emph{linear} data is used as an estimate. Linear data in 1d ($\del_x q_0 = \const$) in \eqref{eq:characteristic} yield for Burgers' equation ($a(q) = q$)
\begin{align}
  \hat \xi  &= \frac{x -  t q_0(x)  + t \del_x q_0(x) x}{1 + t \del_x q_0(x)}  
\end{align}
Usage of this formula as an evolution operator for the pointwise degrees of freedom requires the evaluation of the derivative at a location where the data are not differentiable. Also, equations with more complicated wave speeds lead to a lot more complicated formulae. 

\subsection{Modification of the fixpoint iteration in order to account for shocks} \label{ssec:fixpointmodification}

It is well-known that nonlinear hyperbolic equations develop shocks even when the initial data are smooth. Therefore, when studying the time evolution of the reconstruction, the assumption that no shocks appear cannot always be true. However, the reconstruction is continuous. An initial value problem with continuous data does not develop a shock immediately. The shock can only appear only after a time $t_\text s > 0$. Whenever the time step happens to be small enough ($\Delta t < t_\text s$), the reconstruction did not have time to develop a shock and \eqref{eq:evoscalarsimple} is a good estimate. 

\begin{figure}[h]
 \centering
 \includegraphics[width=\textwidth]{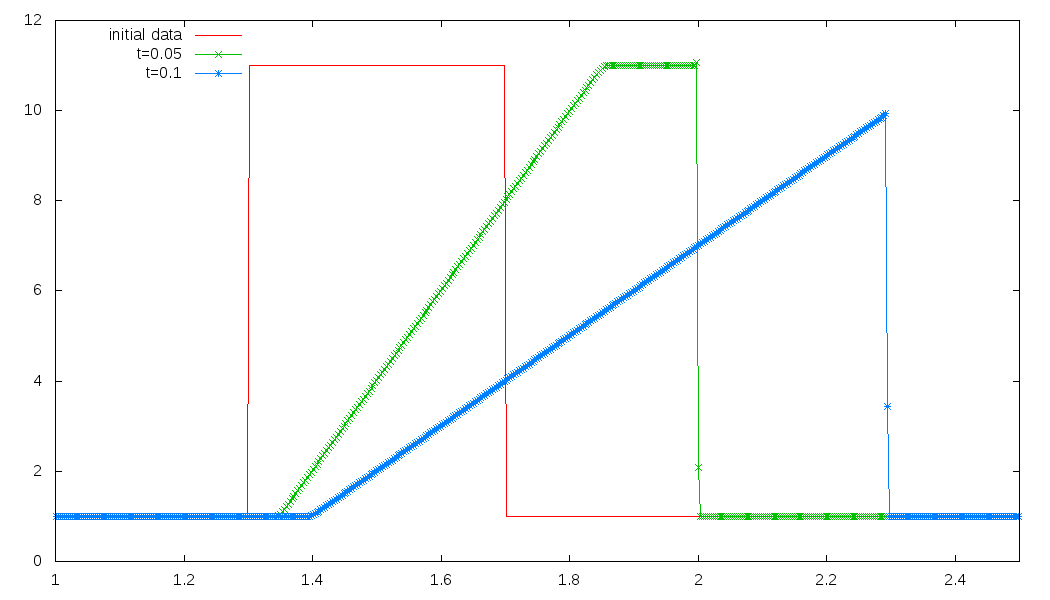}
 \caption{Two Riemann problems for Burgers' equation solved with the Active Flux scheme using iteration \eqref{eq:scalariterationsimple1}--\eqref{eq:scalariterationsimple2} with $\Delta x = 3 \cdot 10^{-3}$. Power law limiting (section \ref{sec:limiting}) has been used. Cell averages are shown.}
 \label{fig:rp-pos}
\end{figure}

This gives an explanation why in certain cases even Riemann problems can be successfully computed with the Active Flux scheme endowed with \eqref{eq:evoscalarsimple}. Fig. \ref{fig:rp-pos} shows such a successful computation of a Riemann problem between values $q_\text{high} = 11$ and $q_\text{low} = 1$ for Burgers' equation. Recall that the initial data in the cell containing the discontinuity are still reconstructed continuously. One can estimate its self-steepening time in this situation as $\frac{q_\text{high} - q_\text{low}}{\Delta x}$. The CFL condition involves the maximum speed $q_\text{high}$ in this case. Thus, for a Riemann problem with uniformly positive values the time step is always smaller than the estimate of the self-steepening time.

Riemann problems involving both positive and negative values do show artefacts. In \cite{kerkmann18} it has been shown, that on such Riemann problems for Burgers' equation evolution operators like the one from \cite{roe17a} fail. \eqref{eq:evoscalarsimple} suffers from very similar problems. In \cite{kerkmann18} it is suggested to revert to a discontinuous reconstruction in this case. However, the failure can be explained by an insufficiently accurate evolution operator rather than tracing it back to continuity of the reconstruction. In order to do this, consider an even simpler Riemann problem for Burgers' equation:
\begin{align}
 q_0(x) = \begin{cases} 1 & x < 0\\ 0 & x > 0 \end{cases} \label{eq:exmaplerpsimplefailure}
\end{align}
The exact solution is a shock moving at speed $\frac12$. However, the evolution operator using \eqref{eq:evoscalarsimple} leaves these data stationary! Indeed, for $x > 0$ the fixpoint iteration is initialized with zero speed. For $x < 0$ the fixpoint iteration converges after one iteration to $\hat \xi^{(n)} = \hat \xi^{(1)} = x - t$ and $q_0(x-t) = 1 \, \forall t$.

Additionally, transonic rarefactions show non-entropic artefacts similar to the ones observed for finite volume schemes with Riemann solvers. Examples of such are shown in Figure \ref{fig:burgersallriemann1char}.

\newcommand{\argmax}{\mathrm{argmax}}

All these problems are removed by modifying the initialization of the fixpoint iteration \eqref{eq:evoscalarsimple} as follows:
\begin{align}
  \hat\vecxi_\ell^{(0)} &= \vec x + \vecdelta_\ell \qquad \ell = 1, \ldots, 2d\label{eq:scalariterationnew1}\\ 
  \hat\vecxi_\ell^{(n)} &= \vec x - \vec a(q_0(\hat \vecxi_\ell^{(n-1)})) t \qquad n = 1, 2, \ldots \label{eq:scalariterationnew2}
 \end{align}
 On two-dimensional Cartesian grids,
 \begin{align}
  \vecdelta_1 &:= \vecc{\Delta x}{0} & \vecdelta_2 &:= \vecc{-\Delta x}{0} & \vecdelta_3 &:= \vecc{0}{\Delta y} & \vecdelta_4 &:= \vecc{0}{-\Delta y}
 \end{align}
 \begin{align}
  L := \argmax_\ell |\vec a(q_0(\hat\vecxi^{(1)}_\ell))|
 \end{align}
and 
\begin{align}
q_\vec x^{n+1} = q_0({\hat\vecxi}^{(2)}_L ) \label{eq:selectingchar}
\end{align}

The reasoning behind this algorithm is the following: Shock formation occurs because of crossing characteristics. Iteration \eqref{eq:scalariterationsimple1}--\eqref{eq:scalariterationsimple2} converges to both footpoints $\hat\vecxi_+$ and $\hat\vecxi_-$ if its initial estimates $\hat\vecxi_\pm^{(0)}$ are chosen appropriately. The above algorithm \eqref{eq:scalariterationnew1}--\eqref{eq:selectingchar} initializes the iteration with the two locations $x \pm \Delta x$, placed symmetrically around $x$. In order to find the solution, one thus needs to estimate which of the characteristics will have survived until time $t$ (and not gone into the shock). The choice of \eqref{eq:selectingchar} is to use the value transported by the quicker characteristic. This choice is inspired by the above example \eqref{eq:exmaplerpsimplefailure} of a Riemann problem, where it makes information flow into the right direction. Of course in general it remains an approximation.

Note that the order of the approximation is not modified, as the modification affects only the initial step of the iteration and is of the order $\mathcal O(\Delta x)$. On smooth solutions, characteristics do not cross, and the two initializations are expected to converge to the same final result. 

Despite its simplicity, in experiments the modification has proven itself able to reliably cure both the artificially stationary shocks and the non-entropic features at transonic rarefactions. One thus may speak of the modification as an entropy fix. To actually prove a statement on the discrete entropy is subject of future work. Instead here a number of different test cases are shown: self-steepening and Riemann problems resulting in strong and weak shocks and (transonic) rarefactions (see sections \ref{ssec:scalarnumerics}--\ref{ssec:otherscalarnumerics}).

\section{Nonlinear systems} \label{sec:systems}
Consider now an $m \times m$ nonlinear hyperbolic system of conservation laws in one spatial dimension:
\begin{align}
 \del_t q + \del_x f(q) &= 0 & q : \mathbb R^+_0 \times \mathbb R &\to \mathbb R^m \label{eq:nonlinsystem1d}\\
 && f : \mathbb R^m \to \mathbb R^m
\end{align}
The Jacobian matrix is denoted by $J(q) := \nabla_q f$. Hyperbolicity guarantees that $J$ has real eigenvalues.

In certain cases a variable change from conservative to characteristic variables $q \mapsto Q$ can be found, such that in the absence of shocks \eqref{eq:nonlinsystem1d} can be rewritten as
\begin{align}
 \del_t Q + \mathrm{diag}(\lambda_1, \lambda_2, \ldots, \lambda_m) \del_x Q &= 0 \label{eq:nonlinsysdiagonal}
\end{align}
with $\lambda_1, \ldots, \lambda_m$ the eigenvalues of $J$.  Denote the initial data as $Q_{i,0}(x) = Q_i(0, x)$, $i=1, \ldots , m$. 

In the linear case this can be solved by solving an advection equation in every component (see e.g. \cite{eymann11}). In the nonlinear case $\lambda_i$ is, in general, a function of all the components of $Q$:
\begin{align}
 \del_t Q_1 + \lambda_1(Q_1, \ldots, Q_m) \del_x Q_1 &= 0\\
 \del_t Q_2 + \lambda_2(Q_1, \ldots, Q_m) \del_x Q_2 &= 0\\
 \nonumber \vdots\\
 \del_t Q_m + \lambda_m(Q_1, \ldots, Q_m) \del_x Q_m &= 0
\end{align}

Therefore, in general the characteristics are curved. This is also a fundamental difference to the nonlinear \emph{scalar} case, where the characteristics remain straight. This is why applying the fixpoint iteration \eqref{eq:evoscalarsimple} to every component of \eqref{eq:nonlinsysdiagonal} does not lead to sufficient order of accuracy (as can be checked by direct computation). A different approximate evolution operator is necessary in the case of systems, which takes into account the curvature of characteristics. Sections \ref{ssec:systemsaveragespeed} and \ref{ssec:rungekutta} describe two such approaches. They yield comparable results, but the strategies and resulting algorithms are fundamentally different. In particular, the algorithm in section \ref{ssec:systemsaveragespeed} does not assume a transformation to characteristic variables. In view of future extensions, e.g. to multiple spatial dimensions, so far it is not clear which of them would be most suitable. They also differ in the nature of necessary computations. Therefore both are presented here.

Even in case the approximate evolution operator is formulated in characteristic variables, the reconstruction still uses conservative variables (as it requires a cell average). At the locations where the initial data need to be evaluated, a transformation to characteristic variables is performed. After obtaining the result of the approximate evolution operator in characteristic variables, they are transformed back to conservative variables.

\subsection{Estimating curved characteristics}
\label{ssec:systemsaveragespeed}

It can be shown by explicit calculation that a straightforward extension of iteration \eqref{eq:scalariterationsimple1}--\eqref{eq:scalariterationsimple2} to \eqref{eq:nonlinsysdiagonal} does not allow to prove a statement analogous to Theorem \ref{thm:fixpointiteration} -- the higher order terms are not correct. This is due to the fact that characteristics are now curved. However, other ways of obtaining a high order estimate can be found. The first is presented in this section, the second -- in section \ref{ssec:rungekutta}

Consider \eqref{eq:nonlinsystem1d} and diagonalize $R J R^{-1} = \Lambda = \mathrm{diag}(\lambda_1, \ldots, \lambda_m)$. Note that, in general, $R$ and all $\lambda_i$ depend on $q$. To make this explicit, write $R(q)$ and $\Lambda(q)$. The equation becomes
\begin{align}
 R(q) \del_t q + \Lambda(q) R(q) \del_x q = 0
\end{align}
Although for some systems it is possible to find $Q(q)$ such that $\del_t Q = R(q) \del_t q$ (as mentioned in the introduction to this section), this is not assumed in the following algorithm.

Denote by $F^{(k)} = R^{-1}\, \mathrm{diag}(0, \ldots, 0,\overset{k}{1},0 \ldots, 0) R$ the projector associated with the $k$-th eigenvalue. Then obviously
\begin{align}
 \sum_{k=1}^m F^{(k)} &= \id & \sum_{k=1}^m F^{(k)} \lambda_k = J
\end{align}

In the following, matrix indices are frequently made explicit, e.g. $q_i$ denotes the $i$-th component of the vector $q$, and $R_{ij}$ the element of $R$ found in its $i$-th row and $j$-th column. It is also sometimes useful to write $\lambda_i(q_1, \ldots, q_m)$ instead of $\lambda_i(q)$.

\begin{theorem}
 Consider a predictor step ($i = 1, \ldots, m$)
\begin{align}
 q^{(i)}_{\beta} := \sum_{k,\alpha = 1}^m F^{(k)}_{\beta \alpha}(x) q_{\alpha,0}\left(x - t \frac{\lambda_i(x) + \lambda_k(x)}{2}\right) 
\end{align}
where $\lambda_i(x)$ is shorthand for $\lambda_i(q_{1,0}(x), \ldots, q_{m,0}(x))$ $\forall i$ and analogously for $F^{(k)}(x)$. Then, with
\begin{align}
 \lambda^*_{i} &:= \lambda_i(q^{(i)}) & R^*_{ij} &:= R_{ij}(q^{(i)}) 
\end{align}
the approximate solution operator 
\begin{align}
 \tilde q_\ell(t, x) := \sum_{i=1}^m (R^*)^{-1}_{\ell i} \sum_{j=1}^m R^*_{ij} q_{j,0}(x - \lambda_i^* t) \label{eq:approxevoopgeneral}
\end{align}
approximates the exact solution of \eqref{eq:nonlinsystem1d} with $J = R^{-1} \mathrm{diag}(\lambda_1, \ldots, \lambda_m) R$ as
\begin{align}q_k(t, x) + \mathcal O(t^3)\end{align}
\end{theorem}
\emph{Note}: The inconspicuously looking equation $R^*_{ij} := R_{ij}(q^{(i)})$ is non-trivial. It states that the rows of $R^*$ are evaluated independently, each on a different predictor value.

\begin{proof}
 Wherever the summation is from $1$ to $m$, it is omitted for the sake of readability. Recall 
 \begin{align} 
 \sum_{i}  R^{-1}_{\ell i} \lambda_i  R_{ij} = J_{\ell j} \label{eq:jacobiandecomp}
 \end{align}
 and note that
 \begin{align}
  \del_t q^{(i)}_{\beta} \Big |_{t=0} &= \sum_{k,\alpha} F^{(k)}_{\beta \alpha}(x) \frac{\lambda_i + \lambda_k}{2} q'_{\alpha,0}(x)\\
  &= -\frac12 \lambda_i  q'_{\beta,0}(x) - \frac12 \sum_{\alpha} J_{\beta \alpha} q'_{\alpha,0}(x) \label{eq:qstarderiv}
 \end{align}
 (The prime denotes differentiation with respect to the unique argument.)
 
 In order to compare the leading order terms in the Taylor series, differentiate the approximate evolution operator with respect to time:
 \begin{align}
 \del_t \tilde q(t, x) &= \sum_{i,j} \del_t \Big ( (R^*)^{-1}_{\ell i}  R^*_{ij} \Big )  q_{j,0}(x - \lambda_i^* t)\\
 &- \sum_{i,j}  (R^*)^{-1}_{\ell i}  R^*_{ij} q'_{j,0}(x - \lambda_i^* t) (\del_t \lambda_i^* t +  \lambda_i^*)
 \end{align}
 On the one hand then
 \begin{align} 
 \del_t \tilde q(t, x) \Big |_{t=0} &= \sum_{i,j} \del_t \Big ( (R^*)^{-1}_{\ell i}  R^*_{ij} \Big ) \Big |_{t=0}  q_{j,0}(x) 
 - \sum_{i,j}  R^{-1}_{\ell i}  R_{ij} q'_{j,0}(x) \lambda_i \\
 &= - \sum_j J_{\ell j} q'_{j,0}(x) 
 \end{align}
 Here $\lambda^*_i |_{t=0} = \lambda_i$, $R^* |_{t=0} = R$ was used. On the other hand
 \begin{align*} 
 \del_t^2 \tilde q(t, x) \Big |_{t=0} &= \sum_{i,j} \del_t^2 \Big ( (R^*)^{-1}_{\ell i}  R^*_{ij} \Big ) \Big |_{t=0} q_{j,0}(x)
  - 2 \sum_{i,j} \del_t \Big ( (R^*)^{-1}_{\ell i}  R^*_{ij} \Big )\Big |_{t=0} q'_{j,0}(x) \lambda_i\\
 &+ \sum_{i,j}  R^{-1}_{\ell i}  R_{ij} q''_{j,0}(x ) \lambda_i^2
 - 2\sum_{i,j}  R^{-1}_{\ell i}  R_{ij} q'_{j,0}(x) \del_t \lambda_i^* \Big |_{t=0} \\ 
 &= \sum_{j} (J^2)_{\ell j}q''_{j,0}(x ) \\
 &-  \sum_j \left( 2\sum_{i} \del_t \Big ( (R^*)^{-1}_{\ell i}  R^*_{ij} \Big )\Big |_{t=0}  \lambda_i 
 + 2\sum_{i,j}  R^{-1}_{\ell i}  R_{ij}  \del_t \lambda_i^* \Big |_{t=0} \right ) q'_{j,0}(x)
 \end{align*}
 The term in brackets can now be expanded using the definitions of $R^*$ and $\lambda^*$:
 \begin{align*}
  & \phantom{m} 2\sum_{i,j} \del_t \Big ( (R^*)^{-1}_{\ell i}  R^*_{ij} \Big )\lambda_i \Big |_{t=0}  
 + 2\sum_{i}  R^{-1}_{\ell i}  R_{ij}  \del_t \lambda_i^* \Big |_{t=0}\\
 &= 2\sum_{i} \del_t  (R^*)^{-1}_{\ell i}\Big |_{t=0}  R_{ij} \lambda_i 
 + 2\sum_{i}   R^{-1}_{\ell i}  \del_t R^*_{ij}\Big |_{t=0}    \lambda_i 
 + 2\sum_{i}  R^{-1}_{\ell i}  R_{ij}  \del_t \lambda_i^* \Big |_{t=0}\\
 \intertext{and using $\del_t R^{-1} = - R^{-1} (\del_t R) R^{-1}$, which follows from $\del_t (R^{-1} R) = \del_t \mathbb{1} = 0$:}
 &=-2 \sum_{i,h,s} R^{-1}_{\ell h} \del_t  R^*_{hs}\Big |_{t=0} R^{-1}_{si}  R_{ij} \lambda_i 
 + 2\sum_{i}   R^{-1}_{\ell i}  \del_t R^*_{ij}\Big |_{t=0}  \lambda_i 
 +2 \sum_{i}  R^{-1}_{\ell i}  R_{ij}  \del_t \lambda_i^* \Big |_{t=0}\\
 \intertext{and with \eqref{eq:qstarderiv}}
 &=\sum_\beta \left( \sum_{i,h,s} R^{-1}_{\ell h} \frac{\del R_{hs}}{\del q_\beta}  \lambda_h  R^{-1}_{si}  R_{ij} \lambda_i 
 - \sum_{i}   R^{-1}_{\ell i} \frac{\del R_{ij}}{\del q_\beta}  \lambda_i^2
 - \sum_{i}  R^{-1}_{\ell i}  R_{ij}  \frac{\del \lambda_i }{\del q_\beta}  \lambda_i \right ) q'_{\beta,0}(x)\\
 &+ \sum_{\alpha,\beta}\left( \sum_{i,h,s} R^{-1}_{\ell h} \frac{\del R_{hs}}{\del q_\beta}   R^{-1}_{si}  \lambda_i R_{ij} 
 - \sum_{i}   R^{-1}_{\ell i}  \lambda_i \frac{\del R_{ij}}{\del q_\beta}   
 - \sum_{i}  R^{-1}_{\ell i}  \frac{\del \lambda_i }{\del q_\beta}  R_{ij}   \right ) J_{\beta \alpha} q'_{\alpha,0}(x)\\
 &=\sum_\beta \left( -\sum_{i,s} J_{\ell s}  \frac{\del R^{-1}_{si}}{\del q_\beta}  \lambda_i  R_{ij}  
 - \sum_{h,s}   J_{\ell s} R^{-1}_{sh}  \lambda_h    \frac{\del R_{hj}}{\del q_\beta} 
 - \sum_{h,s}  J_{\ell s}  R^{-1}_{sh} \frac{\del \lambda_h }{\del q_\beta} R_{hj}    \right ) q'_{\beta,0}(x)\\
 &- \sum_{\alpha,\beta}  \frac{\del J_{\ell j} }{\del q_\beta}     J_{\beta \alpha} q'_{\alpha,0}(x)\\
 &= -\sum_{\beta,s} J_{\ell s} \frac{\del J_{sj}}{\del q_\beta} q'_{\beta,0}(x) - \sum_{\alpha,\beta}  \frac{\del J_{\ell j} }{\del q_\beta}     J_{\beta \alpha} q'_{\alpha,0}(x)
 \end{align*}
 where $\displaystyle \frac{\del R}{\del q_\beta} R^{-1} = - R \frac{\del R^{-1}}{\del q_\beta}$, \eqref{eq:jacobiandecomp} and 
 \begin{align}
 \sum_{h,s} R_{hs} R^{-1}_{si} \lambda_h = \lambda_i
 \end{align}
 was used.
 
 On the other hand, by performing the Cauchy-Kovalevskaya/Lax-Wendroff procedure on the PDE,
  \begin{align*}
   \del_t q_\ell &= - \sum_h J_{\ell h} \del_x q_h\\
   \del_t^2 q_\ell &= -\sum_h \del_t J_{\ell h} \del_x q_h  - \sum_h J_{\ell h} \del_x \del_t q_h\\
   &= - \sum_{h,\beta} \frac{\del J_{\ell h}}{\del q_\beta}  \del_t q_\beta \del_x q_h  + \sum_{h,j} J_{\ell h} \del_x \left(  J_{hj} \del_x q_j \right )\\
   &= \sum_{j,\alpha,\beta} \frac{\del J_{\ell j}}{\del q_\beta} J_{\beta \alpha} \del_x q_\alpha \del_x q_j  + \sum_{h,j,\beta} J_{\ell h} \frac{\del J_{hj}}{\del q_\beta} \del_x q_\beta \del_x q_j  + \sum_j (J^2)_{\ell j} \del_x^2 q_j
  \end{align*}
  which proves the assertion.

\end{proof}\hfill$\Box$

It makes sense to express the matrices $R$ in variables which make the computation simple. In the numerical examples for the full Euler equations, $(\rho, v, p)$ are used:
\begin{align}
 \lambda_+ &= v + c  & \lambda_0 &= v & \lambda_- &= v - c
\end{align}

The transformation matrix $R$ (from $(\rho, v, p)$ to the eigenspace $(+,0,-)$) reads
\begin{align}
 R_{+,\cdot} &= (0,  1, \frac{c }{ \gamma    p} ) \\
 R_{0,\cdot} &= (- \gamma    p  \rho^{- \gamma -1}, 0, \rho^{- \gamma} ) \\
 R_{-,\cdot} &= (0, -1, \frac{c}{ \gamma    p} )
\end{align}

This gives
\begin{align}
 F^{(+)} &= \left( \begin{array}{ccc} 0 &  \frac{ \rho}{ 2c} & \frac{ \rho}{ 2 \gamma  p}  \\ 0 & \frac12 & \frac{c}{2 \gamma  p}  \\ 0 &  \frac{\gamma  p}{2 c} & \frac12 \end{array} \right ) &
 F^{(0)} &= \left( \begin{array}{ccc} 1 & 0 & - \frac{\rho }{ \gamma  p} \\0 & 0 & 0\\ 0 & 0 & 0 \end{array} \right ) & 
 F^{(-)} &= \left( \begin{array}{ccc} 0 & - \frac{\rho }{2c} &  \frac{\rho }{2 \gamma p}  \\ 0 & \frac12 & -\frac{c}{2 \gamma  p}  \\ 0 & - \frac{\gamma p }{2c} & \frac12 \end{array} \right ) 
\end{align}

\begin{lemma}
 When expressing $R$ in the variables $\rho$,$v$ and $p$ for the Euler equations, the approximate evolution operator \eqref{eq:approxevoopgeneral} is exact on contact waves ($p = \const$, $v = \const$).
\end{lemma}
\begin{proof}
 Assume $p = \const$, $v = \const$ uniformly. Then
 
 \begin{align}
 q^{(i)} &= \left( \begin{array}{ccc} 0 &  \frac{ \rho(x)}{ 2c(x)} & \frac{ \rho(x)}{ 2 \gamma  p}  \\ 0 & \frac12 & \frac{c(x)}{2 \gamma  p}  \\ 0 &  \frac{\gamma  p}{2 c(x)} & \frac12 \end{array} \right ) \veccc{\rho \left(x - t \frac{\lambda_i(x) + \lambda_+(x)}{2}\right)}{v}{p}
 \nonumber\\&+ \left( \begin{array}{ccc} 1 & 0 & - \frac{\rho(x) }{ \gamma  p} \\0 & 0 & 0\\ 0 & 0 & 0 \end{array} \right ) \veccc{\rho \left(x - t \frac{\lambda_i(x) + \lambda_0(x)}{2}\right)}{v}{p}
 \nonumber\\&+ \left( \begin{array}{ccc} 0 & - \frac{\rho(x) }{2c(x)} &  \frac{\rho(x) }{2 \gamma p}  \\ 0 & \frac12 & -\frac{c(x)}{2 \gamma  p}  \\ 0 & - \frac{\gamma p }{2c(x)} & \frac12 \end{array} \right )  \veccc{\rho\left(x - t \frac{\lambda_i(x) + \lambda_-(x)}{2}\right)}{v}{p} \\
 &= \veccc{\rho^{(i)}}{v}{p}
\end{align}
with $\rho^{(i)} := \rho \left(x - t \frac{\lambda_i(x) + \lambda_0(x)}{2}\right)$, and the index is taken from the symbolic set $i \in \{ +, -, 0 \}$.

\begin{align}
 \lambda^*_{\pm} &= v \pm c^{(\pm)} & \lambda^*_0 &= v
\end{align}
\begin{align}
 R^* &= \left( \begin{array}{ccc} 0& 1& \frac{c^{(+)} }{ \gamma p} \\ - \gamma    p  (\rho^{(0)})^{- \gamma -1} & 0 & (\rho^{(0)})^{- \gamma} \\ 0 & -1 & \frac{c^{(-)}}{ \gamma    p}  \end{array} \right)\\
 (R^*)^{-1} &= \left( \begin{array}{ccc} \frac{\rho^{(0)}}{c^{(-)} + c^{(+)}} & - \frac{(\rho^{(0)})^{1 + \gamma}}{\gamma p} & \frac{\rho^{(0)}}{c^{(-)} + c^{(+)}} \\ 
 \frac{c^{(-)}}{c^{(-)} + c^{(+)}} & 0 & -\frac{c^{(+)}}{c^{(-)} + c^{(+)}}  \\ 
 \frac{\gamma p}{c^{(-)} + c^{(+)}}   & 0 & \frac{\gamma p}{c^{(-)} + c^{(+)}} \end{array} \right)
\end{align}

\begin{align}
 \tilde q(t, x) &= \veccc{\frac{\rho^{(0)}}{c^{(-)} + c^{(+)}}}{\frac{c^{(-)}}{c^{(-)} + c^{(+)}}}{\frac{\gamma p}{c^{(-)} + c^{(+)}}} (0, 1, \frac{c^{(+)} }{ \gamma p}) \veccc{\rho(x - \lambda_+^* t)}{v}{p}
 \nonumber\\&+\veccc{- \frac{(\rho^{(0)})^{1 + \gamma}}{\gamma p}}{0}{0} (- \gamma    p  (\rho^{(0)})^{- \gamma -1} , 0 , (\rho^{(0)})^{- \gamma}) \veccc{\rho(x - \lambda_0^* t)}{v}{p}
 \nonumber\\&+\veccc{\frac{\rho^{(0)}}{c^{(-)} + c^{(+)}}}{-\frac{c^{(+)}}{c^{(-)} + c^{(+)}} }{\frac{\gamma p}{c^{(-)} + c^{(+)}}} (0 , -1 , \frac{c^{(-)}}{ \gamma    p}) \veccc{\rho(x - \lambda_-^* t)}{v}{p} \\
 &= \veccc{0}{v}{0}  + \veccc{\frac{\rho^{(0)}}{\gamma}}{0}{p}  
 +\veccc{ \rho(x - \lambda_0^* t) - \frac{\rho^{(0)}}{\gamma }  }{0}{0} \\
 &= \veccc{ \rho(x - v t) }{ v}{p} 
\end{align}

This completes the proof.
 
\end{proof}\hfill$\Box$

\begin{corollary}
 For \eqref{eq:nonlinsysdiagonal}, consider a predictor step ($i = 1, \ldots, m$)
\begin{align}
 \hat\xi_{ij}^* &= x - t \frac{\lambda_i(x) + \lambda_j(x)}{2}
\end{align}
where the abbreviation $\lambda_i(x) \equiv \lambda_i(Q_{1,0}(x), \ldots, Q_{m,0}(x))$ has been used.

Then, with
 \begin{align}
 \hat\xi_i &= x - t\lambda_i(Q_{1,0}(\hat\xi_{i1}^*), \ldots, Q_{m,0}(\hat\xi_{im}^*))
 \end{align}
 the approximate solution operator $Q_{i,0}(\hat\xi_i)$ approximates the exact solution as
 \begin{align}Q_{i,0}(\hat\xi_i) = Q_i(t, x) + \mathcal O(t^3)\end{align}
\end{corollary}

Numerical examples are shown in section \ref{ssec:psystem}--\ref{ssec:euler}.

\subsection{Runge-Kutta scheme} \label{ssec:rungekutta}

Recall the second order Runge Kutta method for an ordinary differential equation
\begin{align}
 \dot x = \lambda(t, x)
\end{align}
\begin{align}
 x^* &= x(0) + \alpha t \lambda(0, x(0)) \label{eq:rungekuttasimplestep1}\\
 x(t) &= x(0) + t \left( 1 - \frac{1}{2\alpha} \right ) \lambda(0, x(0)) + t \frac{1}{2\alpha} \lambda(\alpha t, x^*) + \mathcal O(t^3) \label{eq:rungekuttasimplestep2}
\end{align}
For $\alpha = \frac12$ this can be simplified to the midpoint method
\begin{align}
 x(t) = x(0) + t \lambda\left(\frac12 t, x(0) + \frac12 t \lambda(0, x(0))\right) + \mathcal O(t^3)
\end{align}

For simplicity of presentation, consider $m=2$. The Runge-Kutta integration can be applied to the characteristic relations
\begin{align}
 \xi_1' &= \lambda_1(Q_1(t, \xi_1), Q_2(t, \xi_1)) \label{eq:charrelationssystem1}\\
 \xi_2' &= \lambda_2(Q_1(t, \xi_2), Q_2(t, \xi_2)) \label{eq:charrelationssystem2}
\end{align}
that govern the time evolution of the characteristic curves $\xi_i : \mathbb R^+_0 \to \mathbb R$, $i = 1, 2$.

\begin{theorem}
Consider the predictor step
\begin{align}
 \hat\xi_1^* &:= x - \alpha t \lambda_1(Q_{1,0}(x), Q_{2,0}(x)) \label{eq:rungekuttastep1}\\
 \hat\xi_2^* &:= x - \alpha t \lambda_2(Q_{1,0}(x), Q_{2,0}(x)) \label{eq:rungekuttastep12}
\end{align}
Then define
\begin{align}
 \lambda_1^* := \lambda_1\Big (&Q_{1,0}(\hat\xi_1^* - \alpha t \lambda_1(Q_{1,0}(\hat\xi_1^*), Q_{2,0}(\hat\xi_1^*)) ), \\ 
 &  Q_{2,0}(\hat\xi_1^* - \alpha t \lambda_2(Q_{1,0}(\hat\xi_1^*), Q_{2,0}(\hat\xi_1^*))) \Big) \label{eq:rungekuttastep2}\\
 \lambda_2^* := \lambda_2\Big (&Q_{1,0}(\hat\xi_2^* - \alpha t \lambda_1(Q_{1,0}(\hat\xi_2^*), Q_{2,0}(\hat\xi_2^*)) ),  \\
 &  Q_{2,0}(\hat\xi_2^* - \alpha t \lambda_2(Q_{1,0}(\hat\xi_2^*), Q_{2,0}(\hat\xi_2^*))) \Big) \label{eq:rungekuttastep22}
\end{align}
and 
\begin{align}
 \hat\xi_1 &:= x - t \left( 1 - \frac{1}{2\alpha} \right ) \lambda_1(Q_{1,0}(x), Q_{2,0}(x)) - t \frac{1}{2\alpha} \lambda_1^*\\
 \hat\xi_2 &:= x - t \left( 1 - \frac{1}{2\alpha} \right ) \lambda_2(Q_{1,0}(x), Q_{2,0}(x)) - t \frac{1}{2\alpha} \lambda_2^*
\end{align}

Now, $Q_{1,0}(\hat\xi_1), Q_{2,0}(\hat\xi_2)$ approximates $Q_1(t,x), Q_2(t,x)$ with an error $\mathcal O(t^3)$ for any $\alpha \in (0, 1]$.
\end{theorem}

\begin{proof}
 Apply the Runge-Kutta scheme \eqref{eq:rungekuttasimplestep1}--\eqref{eq:rungekuttasimplestep2} to the \emph{backward} time evolution of \eqref{eq:charrelationssystem1}--\eqref{eq:charrelationssystem2} with $\xi_i(t) = x$, $i=1,2$. The equation to solve here is
 \begin{align}
  \vecc{\xi_1'}{\xi_2'} = \vecc{\lambda_1(Q_1(t, \vec \xi_1),Q_2(t, \vec \xi_1) )}{\lambda_2(Q_1(t, \vec \xi_2),Q_2(t, \vec \xi_2) )}
 \end{align}
 The predictor step \eqref{eq:rungekuttasimplestep1} reads
 \begin{align}
  \vecc{\xi_1(t)}{\xi_2(t)} - \vecc{\xi_1(t-\alpha t)}{\xi_2(t-\alpha t)} = \alpha  t \vecc{\lambda_1(Q_1(0, \vec \xi_1(t)),Q_2(0, \vec \xi_1(t)) )}{\lambda_2(Q_1(0, \vec \xi_2(t)),Q_2(0, \vec \xi_2(t)) )}
 \end{align}
 which gives \eqref{eq:rungekuttastep1}--\eqref{eq:rungekuttastep12} defining $\xi_i(t(1-\alpha)) =: \hat\xi^*_i$, $i=1,2$. Now in \eqref{eq:rungekuttasimplestep2} the speeds at time $t(1-\alpha)$ are used. One thus needs an estimate of $Q_i$, $i=1,2$ at time $t(1-\alpha)$. $Q_1$ is constant along the $i$-th characteristic. Thus, for any location $\xi$
 \begin{align}
  Q_i(t(1-\alpha), \xi) = Q_i\Big(0, \xi - \alpha t \lambda_i(Q_1(0, \xi), Q_2(0, \xi))\Big)
 \end{align}
 is an estimate of the solution, which yields \eqref{eq:rungekuttastep2}--\eqref{eq:rungekuttastep22}.
\end{proof}\hfill$\Box$

{\sl Note}: The extension to any $m$ is obtained analogously.

By analogy with the modified fixpoint iteration \eqref{eq:scalariterationnew1}--\eqref{eq:scalariterationnew2} the following modification is suggested for the case of systems: Instead of \eqref{eq:rungekuttastep1}--\eqref{eq:rungekuttastep12}, compute ($\ell = 1, 2$)
\begin{align}
 \hat\xi_{1,\ell}^* &:= x  - \alpha t \lambda_1(Q_{1,0}(x + \delta_\ell), Q_{2,0}(x+ \delta_\ell)) \label{eq:rungekuttastep1modif}\\
 \hat\xi_{2,\ell}^* &:= x  - \alpha t \lambda_2(Q_{1,0}(x+ \delta_\ell), Q_{2,0}(x+ \delta_\ell)) \label{eq:rungekuttastep12modif}
\end{align}
\begin{align}
 \delta_1 &= \Delta x & \delta_2 &= -\Delta x
\end{align}
For each $\ell$, proceed with the algorithm to obtain $\hat \xi_{1,\ell}, \hat \xi_{2,\ell}$. Choose
\begin{align}
 \hat \xi_1 &= \begin{cases} \hat \xi_{1,1} & \text{if } |\hat \xi_{1,1} - x| > |\hat \xi_{1,2} - x|\\ \hat \xi_{1,2} & \text{else} \end{cases} & 
 \hat \xi_2 &= \begin{cases} \hat \xi_{2,1} & \text{if }|\hat \xi_{2,1} - x| > |\hat \xi_{2,2} - x|\\ \hat \xi_{2,2} & \text{else} \end{cases} \label{eq:rungekuttamodifselect}
\end{align}
Define the approximate solution operator by $(Q_{1,0}(\hat \xi_1), Q_{2,0}(\hat \xi_2))$.

Numerical examples are shown in section \ref{ssec:psystem}--\ref{ssec:euler}.

\section{Limiting in one spatial dimension}\label{sec:limiting}

The Active Flux scheme uses a conservative reconstruction in order to evolve the pointwise degrees of freedom. It has to fulfill several conditions (i.e. \eqref{eq:interpolpropertyrecon}--\eqref{eq:conservpropertyrecon}). In one spatial dimension conservation and interpolation of the two point values at cell boundaries amount to three conditions. The natural choice therefore is a parabola (e.g. in \cite{vanleer77}). However, polynomials in general do not fulfill a maximum principle: The maximum of the reconstruction $q_\text{recon}(x)$ in cell $i$ can exceed $\max(\bar q_i, q_{i+\frac12}, q_{i-\frac12})$. To correct this (whenever possible) is the objective of the limiting procedure introduced below.

The starting point is an analysis of the failure of the parabolic reconstruction to be monotone. Assume in the following that $q_{i-\frac12} < q_{i+\frac12}$; for the opposite situation analogous statements are true.
\begin{theorem}\label{thm:paraboliccases}
 \begin{enumerate}[i)]
 \item If $q_{i-\frac12} < \bar q_i < q_{i+\frac12}$ then there exists a monotone continuous function satisfying \eqref{eq:interpolpropertyrecon}--\eqref{eq:conservpropertyrecon}.
 \item With $r := \frac{q_{i+\frac12} -q_{i-\frac12} }{3} > 0$, if $\bar q_i > q_{i+\frac12} - r$ or $\bar q_i < q_{i-\frac12} + r$ then the parabolic reconstruction \eqref{eq:parabolicrecon} is not monotone. \label{it:paraboliccases2}
 \end{enumerate}
 \end{theorem} 
\begin{proof}
 \begin{enumerate}[i)]
  \item E.g. a piecewise linear function can easily be constructed to fulfill the conditions.
  \item The reconstruction \eqref{eq:parabolicrecon} has a maximum inside $[-\frac{\Delta x}{2}, \frac{\Delta x}{2}]$ if the average is too close to the point values. The maximum is located at $\frac{\Delta x (q_{i+\frac12} - q_{i-\frac12})}{12 \bar q_i - 6 (q_{i-\frac12} + q_{i+\frac12})}$. Equating this to $\pm\frac{\Delta x}{2}$ yields the bounds.
 \end{enumerate}
\end{proof}\hfill$\Box$

\begin{figure}[h]
 \centering
 \includegraphics[width=0.75\textwidth]{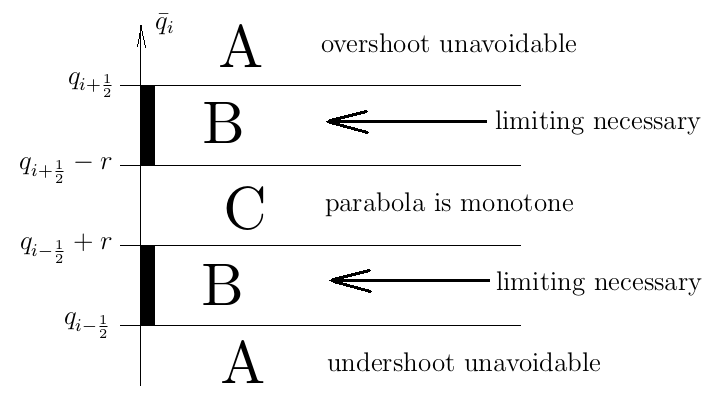}
 \caption{Values of $\bar q_i$ for which limiting is possible and necessary.}
 \label{fig:limitingsketch}
\end{figure}

Thus, there are three possible cases, which are shown in Figure \ref{fig:limitingsketch}.
\begin{enumerate}[A.]
 \item $\bar q_i > q_{i+\frac12}$ or $\bar q_i < q_{i-\frac12}$: no continuous monotone reconstruction exists: use the parabolic reconstruction.
 \item $q_{i-\frac12} < \bar q_i < q_{i-\frac12} + r$ or $q_{i+\frac12}-r < \bar q_i < q_{i+\frac12}$: correction is needed.
 \item $q_{i-\frac12} + r \leq \bar q_i \leq q_{i+\frac12} - r$: the parabolic reconstruction is monotone and no limiting needed.
\end{enumerate} 
 
The following function corrects the failure of the parabolic reconstruction:

\begin{theorem}[Power law limiting] \label{thm:powerlawlimiting}
 Under the conditions of Theorem \ref{thm:paraboliccases}\ref{it:paraboliccases2}) the function
 \begin{align}
    p_N(x) &= q_{i-\frac12} + (q_{i+\frac12} - q_{i-\frac12}) \left( \frac{x - x_i + \Delta x / 2}{\Delta x}   \right )^N  \qquad x_i-\frac{\Delta x}{2} < x < x_i + \frac{\Delta x}{2}  \label{eq:powerlaw}
 \end{align}
 with $\displaystyle N = \frac{q_{i+\frac12} - \bar q_i}{\bar q_i - q_{i-\frac12}}$ fulfills \eqref{eq:interpolpropertyrecon}--\eqref{eq:conservpropertyrecon} and is monotone.
\end{theorem}
\begin{proof}
 Monotonicity and \eqref{eq:interpolpropertyrecon} are obvious. For \eqref{eq:conservpropertyrecon} compute
 \begin{align}
  \frac{1}{\Delta x} \int_{x_i-\frac{\Delta x}{2}}^{x_i + \frac{\Delta x}{2}} \dd x \, p_N(x) &= q_{i-\frac12} + 
  (q_{i+\frac12} - q_{i-\frac12}) \frac{1}{\Delta x^{N+1}} \int_{0}^{\Delta x} \dd x \,x^N\\
  &= q_{i-\frac12} + 
  (q_{i+\frac12} - q_{i-\frac12}) \frac{1}{N+1} = \bar q_i
 \end{align}
\end{proof}\hfill$\Box$

\begin{figure}[h]
 \centering
 \includegraphics[width=0.32\textwidth]{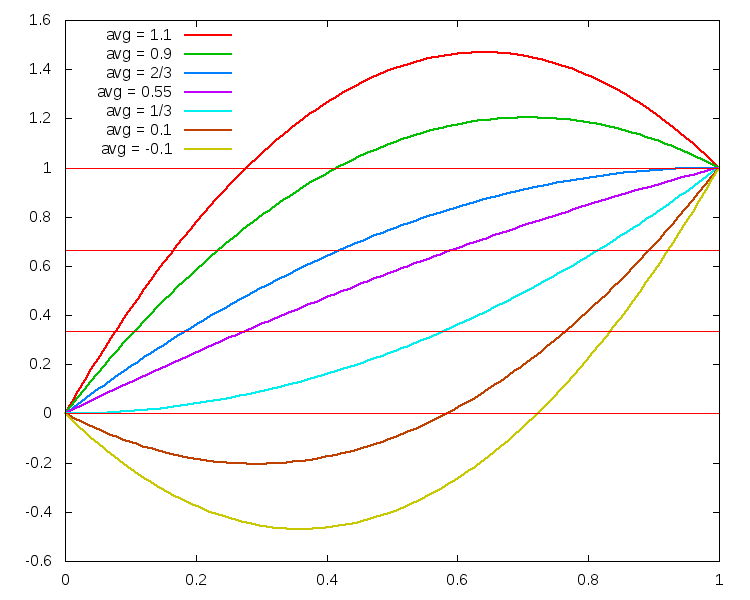} \hfill \includegraphics[width=0.32\textwidth]{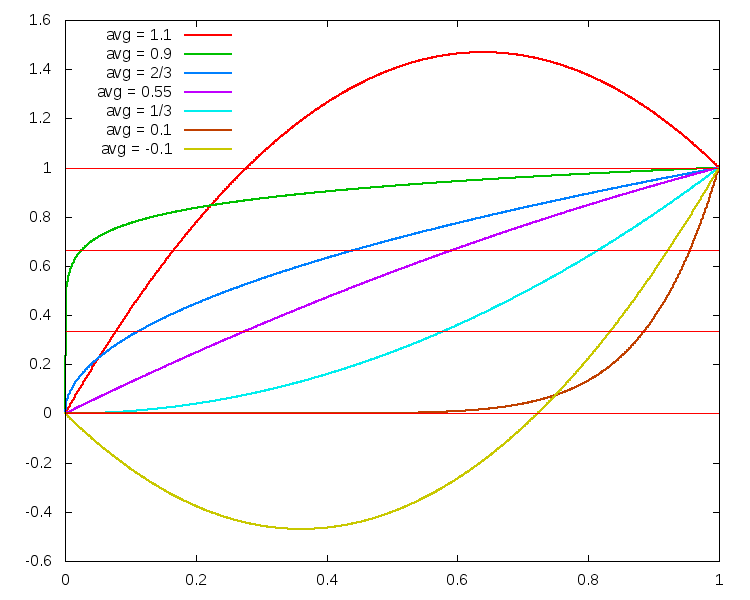} \hfill \includegraphics[width=0.32\textwidth]{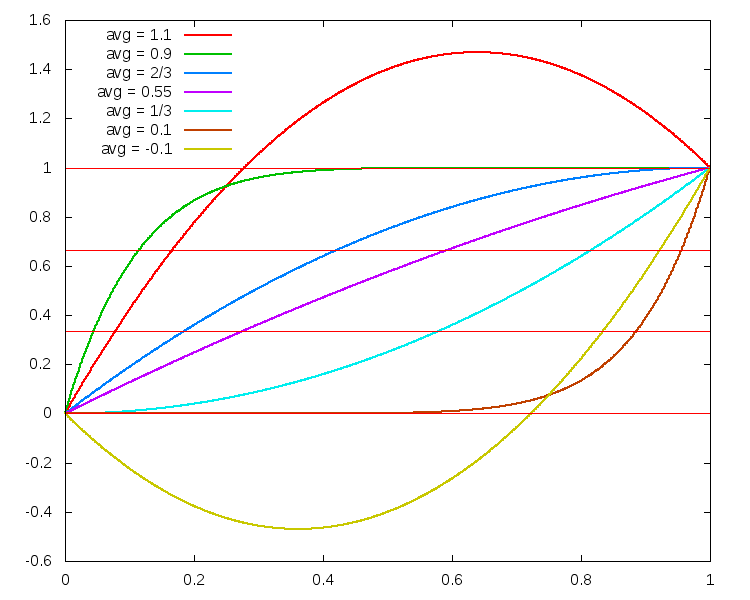}
 \caption{Example of reconstructions with $\Delta x = 1$, $x_i = 0.5$. The data are $q_{i-\frac12} = 0$, $q_{i+\frac12} = 1$ and different averages: $\bar q_i \in \{ 1.1, 0.9, \frac23, 0.55, \frac13, 0.1, -0.1 \}$. {\sl Left}: Parabolic reconstruction \eqref{eq:parabolicrecon}. It is only monotone for $\frac13 \leq \bar q_i \leq \frac23$. {\sl Center}: Limiting applied. Now the reconstruction is monotone for $q_{i-\frac12} \leq \bar q_i \leq q_{i+\frac12}$. Outside this range no monotone function can be found, and parabolic reconstruction is still used. {\sl Right}: A version of the limiting symmetric with respect to $N \mapsto \frac{1}{N}$ according to \eqref{eq:limitingsymmetrizedpowerlaw}.}
 \label{fig:limiting}
\end{figure}

Thus, using the limiter amounts to replacing the parabolic reconstruction by \eqref{eq:powerlaw} in the region B:
\begin{align}
 q_\text{recon}(x) = \begin{cases} p_N(x) & (q_{i-\frac12} < \bar q_i < q_{i-\frac12} + r) \text{ or }(q_{i+\frac12} - r < \bar q_i < q_{i+\frac12})\\
                       \text{parabolic \eqref{eq:parabolicrecon}} & \text{else} \label{eq:limitingonepowerlaw}
                     \end{cases}
\end{align}
The effect is illustrated in Figure \ref{fig:limiting}.

{\sl Notes:} 
\begin{enumerate}[i)]
\item If $\bar q_i > q_{i-\frac12} + r$ then $N < 2$. Therefore inside the region C ($q_{i-\frac12} + r < \bar q_i < q_{i+\frac12} - r)$ where the parabolic reconstruction \eqref{eq:parabolicrecon} is monotone, the power law would have a formal approximation order less than the parabola.
\item Usage of the parabolic reconstruction whenever an overshoot (undershoot) is unavoidable (region A) is making the limiter affect maxima (minima) as little as possible, and thus to avoid clipping. In practice, limiting is discarded if it would imply $\max(N, \frac{1}{N}) > 50$.
\item As can be seen from Figure \ref{fig:limiting} (center), the curves for $N$ and $1/N$ are not symmetric. One thus can consider
\begin{align}
 q_\text{recon}(x) = \begin{cases} p_N(x) & \text{if }q_{i-\frac12} < \bar q_i < q_{i-\frac12} + r\\
                       q_{i+\frac12} - (q_{i+\frac12} - q_{i-\frac12}) \left( \frac{\Delta x / 2 - x + x_i}{\Delta x}   \right )^{1/N }
                       & \text{if }q_{i+\frac12} - r < \bar q_i < q_{i+\frac12}\\
                       \text{parabolic \eqref{eq:parabolicrecon}} & \text{else} 
                     \end{cases} \label{eq:limitingsymmetrizedpowerlaw}
\end{align}
as a limiting strategy instead of \eqref{eq:limitingonepowerlaw}. The result is shown in Fig. \ref{fig:limiting} (right).
\item For $q_{i-\frac12} > q_{i+\frac12}$ the parabola is monotone if $\bar q_i$ fulfills
\begin{align}
 q_{i+\frac12} - \frac{q_{i-\frac12} - q_{i+\frac12}}{3} \leq \bar q_i \leq q_{i-\frac12} - \frac{q_{i-\frac12} - q_{i+\frac12}}{3}
\end{align} 
but the formula \eqref{eq:powerlaw} remains unchanged.

\end{enumerate}

\begin{figure}[h]
 \centering
 \includegraphics[width=\textwidth]{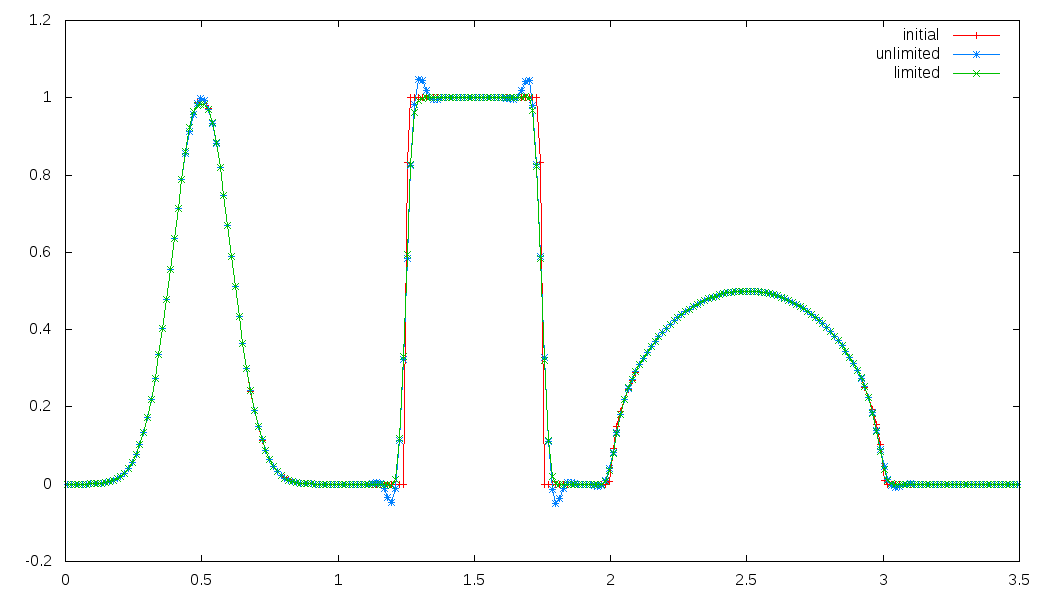}
 \caption{Effect of limiting for linear advection with speed 1.0. Here, $\Delta x = 3.5/250$, the CFL number is chosen $0.45$ and the curves show cell averages of the initial data and of the solution at $t=3.5$, i.e. after one revolution with periodic boundaries. The blue curve shows the solution without limiting, and the green one with limiting.}
 \label{fig:limitingadvection}
\end{figure}

In Figure \ref{fig:limitingadvection} the effect of limiting is shown for linear advection. The initial data are compared to the numerical solution on a periodic grid after one revolution. 

The approximate evolution operators in sections \ref{sec:scalar} and \ref{sec:systems} approximate the solution at time $t$ by evaluating the initial data at some particular location. The above limiting procedure thus guarantees that the point value update satisfies the maximum principle in all cases when a monotone reconstruction is at all possible. This does not mean, however, that the full numerical solution will be free of spurious oscillations. The finite volume update of the average might still lead to the appearance of oscillations. In practice, for linear advection one does observe oscillations for large CFL numbers, but they are much smaller than without limiting. A limiting of the finite volume step for Active Flux has not yet been considered in the literature and is subject of future work.

For systems, limiting is applied to conservative quantities individually. Numerical examples of limiting for nonlinear problems are shown in section \ref{sec:numerical}.

\section{Numerical examples} \label{sec:numerical}

The Active Flux scheme endowed with the approximate solution operators of the above section is now applied to several problems in order to assess its abilities experimentally. First, scalar equations and systems in one spatial dimension are considered (sections \ref{ssec:scalarnumerics}--\ref{ssec:euler}); in section \ref{ssec:scalarnumericsmultid} the scheme is applied to multi-dimensional scalar conservation laws.

The CFL condition is applied using the maximum absolute value of the characteristic speed (eigenvalue of the Jacobian in case of systems) evaluated at the point values (edge midpoint values in the multi-dimensional case).

Third order accuracy in time requires the computation of the point values at both the full and the half time step. Simpson's rule in time then is used to compute the fluxes necessary for the cell averages update \eqref{eq:finitevolume}.

\subsection{Burgers' equation} \label{ssec:scalarnumerics}

Here, Burgers' equation
\begin{align}
 \del_t q + \del_x \left(  \frac{q^2}{2}  \right ) &= 0
\end{align}
is solved with the Active Flux scheme. In all cases the approximate evolution operator \eqref{eq:selectingchar} is used with $a(q) = q$. Additionally, section \ref{ssec:burgersriemannproblems} shows the artefacts appearing upon usage of the simple, unmodified fixpoint iteration \eqref{eq:evoscalarsimple}, and their absence when using \eqref{eq:selectingchar}.

\subsubsection{Convergence study}

\begin{figure}[h]
 \centering
 \includegraphics[width=0.48\textwidth]{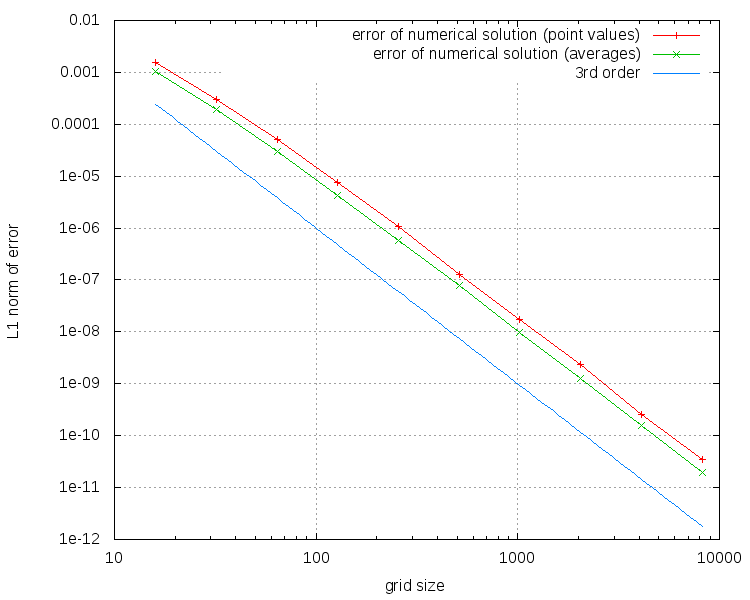}\hfill\includegraphics[width=0.48\textwidth]{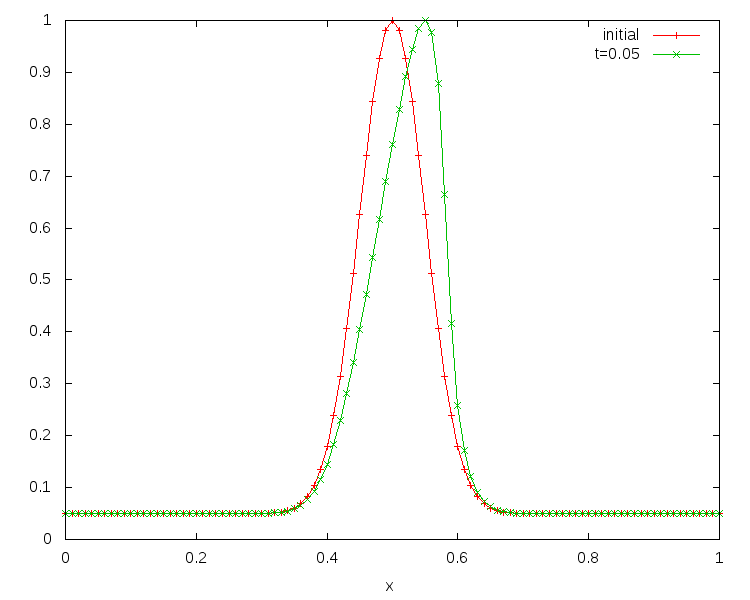}
 \caption{Burgers' equation solved with Active Flux using the modified fixpoint iteration \eqref{eq:selectingchar}. \textit{Left}: Third order convergence of the numerical solution on both point values and averages. \textit{Right}: Setup and numerical solution for $\Delta x = 1/100$. The error is evaluated at $t=0.05$. The CFL number is 0.45.}
 \label{fig:burgersconvergence}
\end{figure}

In order to assess the experimental convergence rate, Gaussian initial data are evolved on grids of different refinement. The simulation is stopped before a shock occurs. Figure \ref{fig:burgersconvergence} (right) shows the setup and the solution at final time on a grid with $\Delta x = 1 /100$. Figure \ref{fig:burgersconvergence} (left) shows the $\ell_1$ norm of the error of both the point values and the averages at $t=0.05$. The reference solution is obtained by evolving a piecewise linear approximation of the data exactly on a grid 30 times finer. For the error of the averages, Simpson's rule is used to compute cell averages of the reference solution. Periodic boundaries are used, though the Gaussian decays sufficiently quickly to a constant towards the boundaries. Limiting is not used.

\newpage

\subsubsection{Self-steepening}

\begin{figure}[h]
 \centering
 \includegraphics[width=0.4\textwidth]{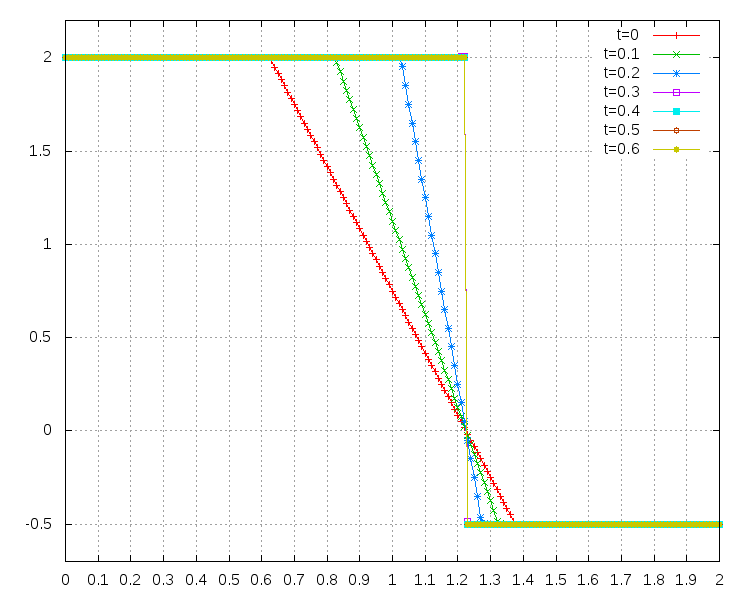} \hfill \includegraphics[width=0.4\textwidth]{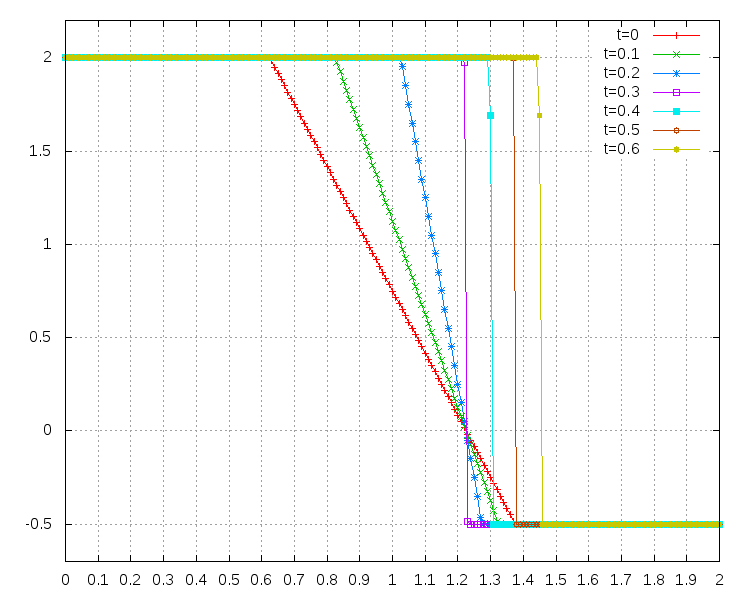}
 \caption{Burgers' equation evolved with the Active Flux scheme. The initial data self-steepen and a shock forms at $t=0.3$. Point values are shown at times $t \in \{ 0, 0.1, \ldots, 0.6 \}$. The solution is computed on a much larger grid such that the boundaries are of no influence. $\Delta x = 1/100$ and the CFL number is 0.9. Limiting procedure according to Theorem \ref{thm:powerlawlimiting} / Equation \eqref{eq:limitingonepowerlaw} is used. \textit{Left}: Usage of fixpoint iteration \eqref{eq:evoscalarsimple} yields zero shock speed (the corresponding lines are on top of each other). \textit{Right}: The modified iteration \eqref{eq:selectingchar} yields the correct shock speed.}
 \label{fig:burgersinverserf} 
\end{figure}

\begin{figure}[h]
 \centering
 \includegraphics[width=0.68\textwidth]{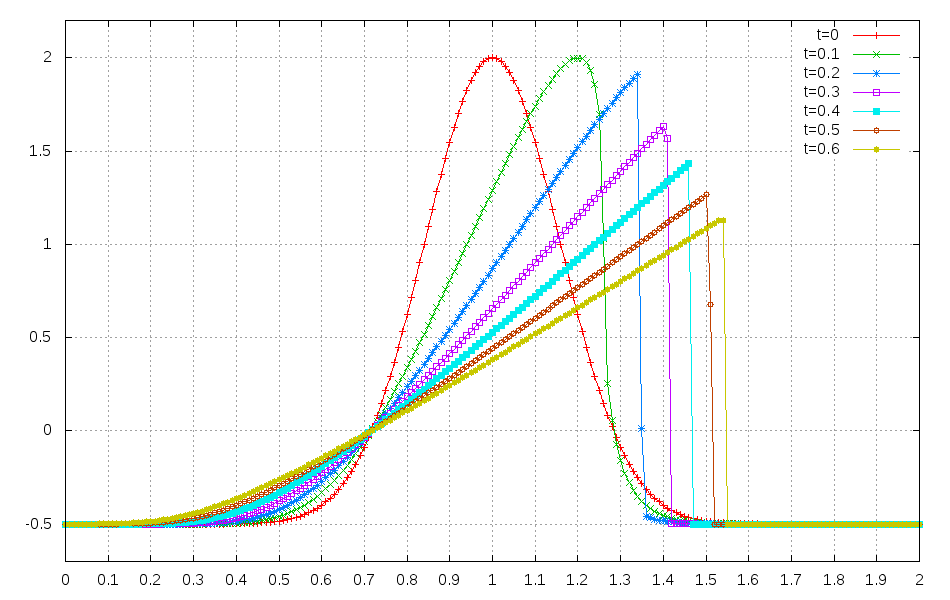}
 \caption{Burgers' equation evolved with the Active Flux scheme. Point values are shown at times $t \in \{ 0, 0.1, \ldots, 0.6 \}$. The solution is computed on a much larger grid such that the boundaries are of no influence. $\Delta x = 1/100$ and the CFL number is 0.9. Limiting procedure according to Theorem \ref{thm:powerlawlimiting} / Equation \eqref{eq:limitingonepowerlaw} is used.}
 \label{fig:burgersgaussianselfsteep} 
\end{figure}

\begin{figure}[h]
 \centering
 \includegraphics[width=0.48\textwidth]{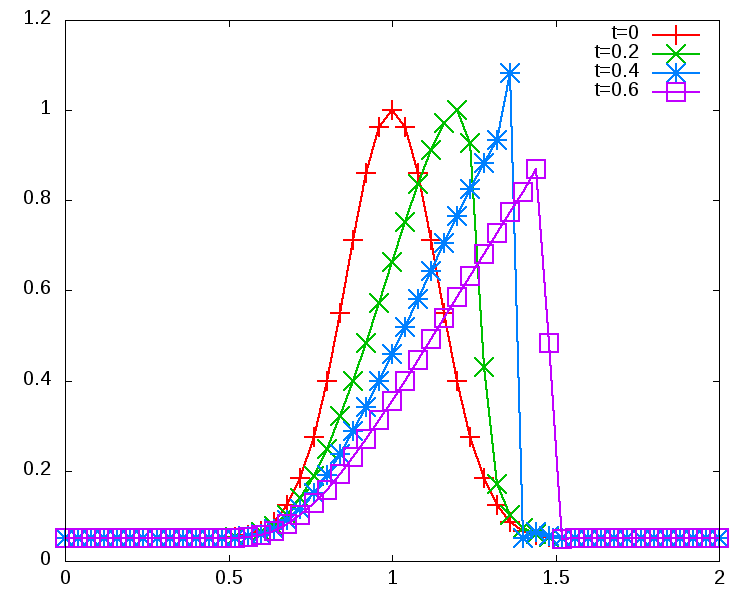} \hfill \includegraphics[width=0.48\textwidth]{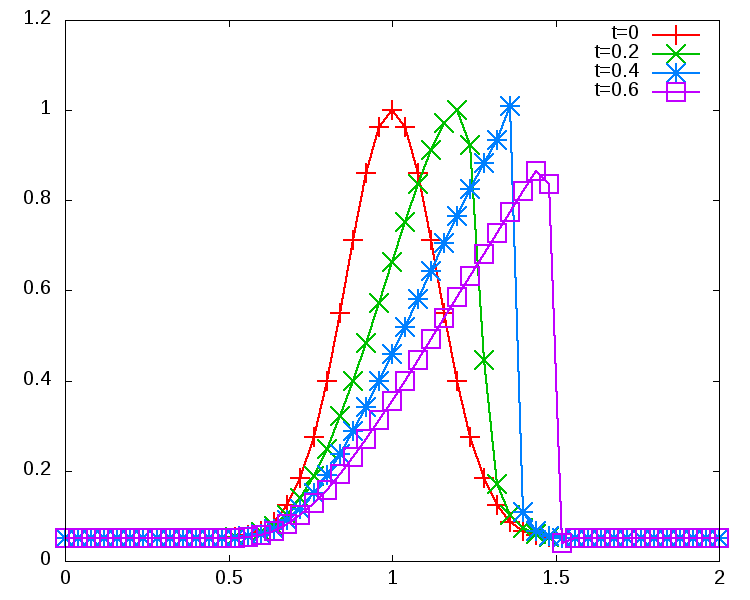} 
 \caption{Burgers' equation evolved with the Active Flux scheme. Point values are shown at times $t \in \{ 0, 0.2, 0.4, 0.6 \}$. $\Delta x = 1/25$ and the CFL number is 0.45. No limiting is used. \emph{Left}: Simple fixpoint iteration \eqref{eq:evoscalarsimple} is producing a spurious oscillation. \emph{Right}: Modified fixpoint iteration \eqref{eq:selectingchar} yields much better results, even without limiting.}
 \label{fig:burgersgaussianosci} 
\end{figure}

Figures \ref{fig:burgersinverserf} and \ref{fig:burgersgaussianselfsteep} show continuous initial data which self-steepen into a shock. In both cases the arising shock connects a positive and a negative value -- a situation in which the simple fixpoint iteration \eqref{eq:evoscalarsimple} is known to fail. Figure \ref{fig:burgersinverserf} (left) shows such a failure: the shock is stationary, which violates the Rankine-Hugoniot condition (this problem has been first reported in \cite{kerkmann18}). Figure \ref{fig:burgersinverserf} (right) demonstrates that the modified iteration \eqref{eq:selectingchar} yields the correct shock speed $0.75$. Moreover, in figure \ref{fig:burgersgaussianosci} it is shown that the simple fixpoint iteration \eqref{eq:evoscalarsimple} produces a spurious oscillation when a shock appears. This is not the usual oscillation due to the high order of the method, as it is removed upon usage of the modified fixpoint iteration \eqref{eq:selectingchar}, even if limiting is not applied. The reason for the appearance of the oscillation is related to the stationary-shock artefact of figure \ref{fig:burgersinverserf} (left). As is emphasized in \cite{kerkmann18}, stationary point values imply fluxes which do not change in time. At the location of the shock, the two fluxes of the cell are different and the average keeps growing. This pile-up is observed as an oscillation in figure \ref{fig:burgersgaussianosci}, although the non-zero value on the left side of the shock is enough to make it travel at the right speed.

\subsubsection{Riemann problems} \label{ssec:burgersriemannproblems}

\begin{figure}[h]
 \centering
 \includegraphics[width=\textwidth]{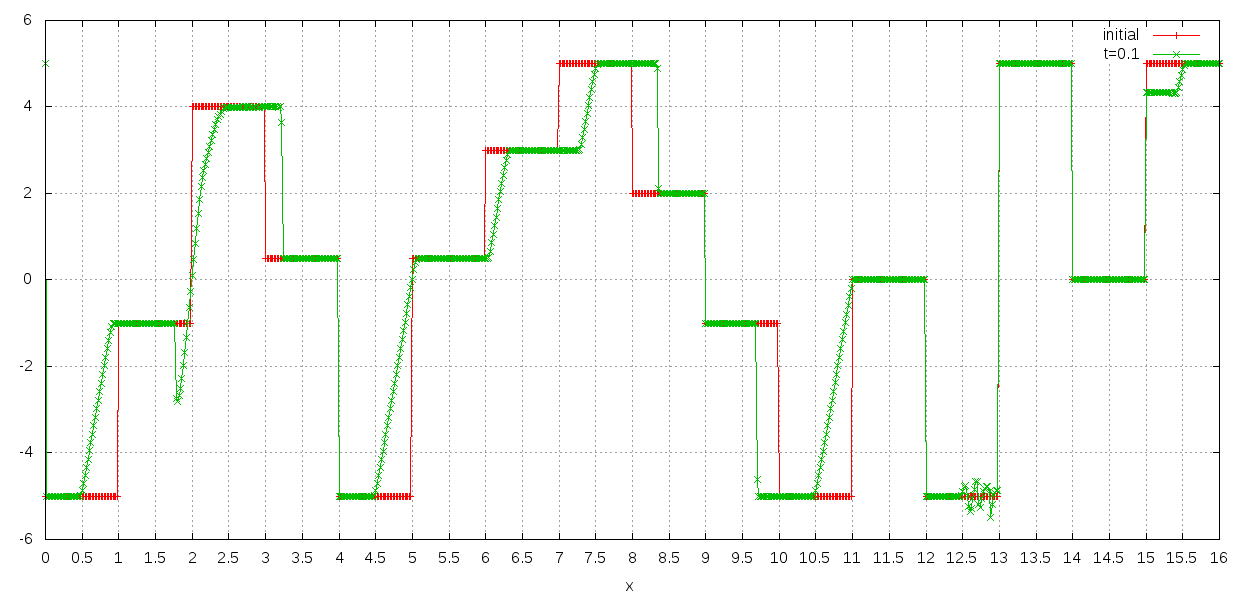}
 \caption{Burgers' equation evolved with the Active Flux scheme using the simple fixpoint iteration \eqref{eq:evoscalarsimple}: wrong shock speeds and artefacts in transonic rarefactions are visible. The initial data are piecewise constant. Point values of the solution are shown at time $t = 0.1$. The boundaries are periodic. $\Delta x = 2/100$ and the CFL number is 0.45. Limiting procedure according to Theorem \ref{thm:powerlawlimiting} / Equation \eqref{eq:limitingonepowerlaw} is used.}
 \label{fig:burgersallriemann1char} %
\end{figure}

Finally, a number of Riemann problems are used in order to assess the properties of the suggested modification of the fixpoint iteration. Fig. \ref{fig:burgersallriemann1char} shows an initial setup chosen to include many interesting cases at once, including strong and weak shocks and rarefactions, and shocks and rarefactions where the discontinuity crosses 0. These latter are known to be the most problematic. On the right also some of the constant states are 0. The boundaries are periodic. 

In the same Figure, a solution at time $t=0.1$ is shown which is using the simple, unmodified fixpoint iteration \eqref{eq:evoscalarsimple}. One observes that only uniformly positive or uniformly negative shocks are evolved correctly. Some of the rarefactions are correct, some have nonentropic features.

\begin{figure}[h]
 \centering
 \includegraphics[width=\textwidth]{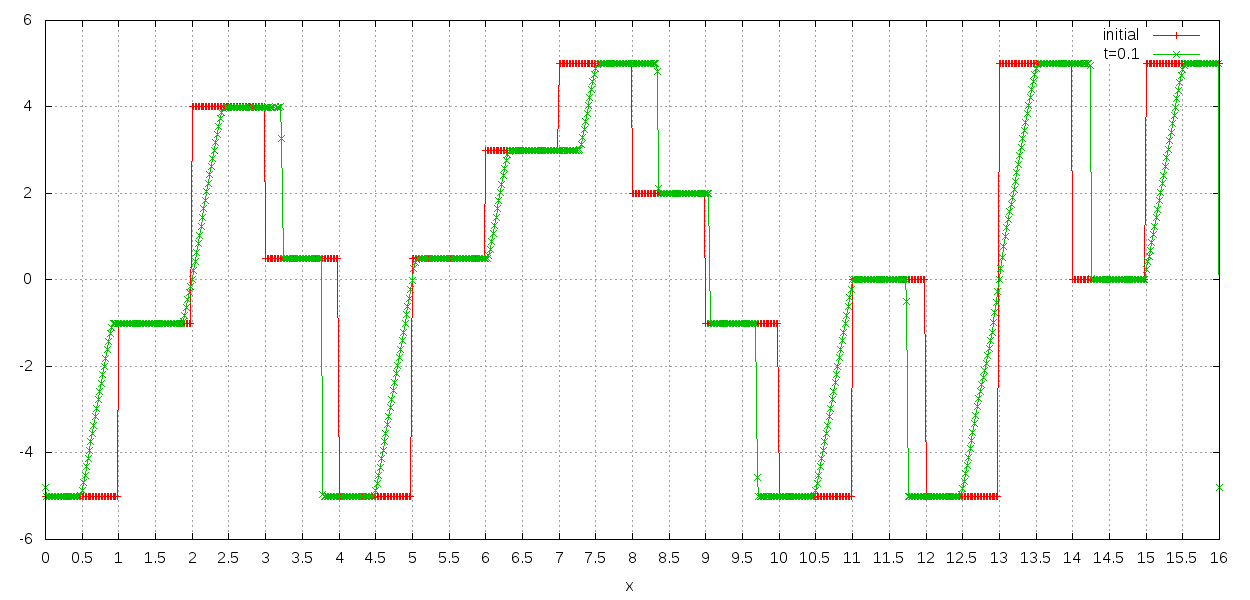}
 \caption{Burgers' equation evolved with the Active Flux scheme using the modified fixpoint iteration \eqref{eq:selectingchar}. All the shock speeds are correct and the rarefactions are not distorted by non-entropic shocks. Point values of the solution are shown at time $t = 0.1$. Setup as in Figure \ref{fig:burgersallriemann1char}.}
 \label{fig:burgersallriemann} %
\end{figure}

For comparison, Fig. \ref{fig:burgersallriemann} shows the evolution of the same setup using the modified fixpoint iteration \eqref{eq:selectingchar}. Now all the shocks have the correct speed, and the rarefactions do not contain additional non-entropic shocks. The strong rarefaction with values symmetric around 0 seems particularly difficult to capture. For large CFL numbers, little artefacts have been observed, which, however, vanish upon grid refinement. 

To the author's knowledge the performance of Active Flux for nonlinear problems has never been studied on transonic rarefactions or strong shocks in the literature.

\clearpage

\subsection{Other scalar equations} \label{ssec:otherscalarnumerics}

Consider the following scalar conservation law
\begin{align}
 \del_t q + \del_x \left( \frac{q^4}{4} \right ) &= 0 \label{eq:quarticequation}
\end{align}

The speed of a shock joining $q_\text L$ and $q_\text R$ is
\begin{align}
 s = \frac{q_\text L^3 +q_\text R q_\text L^2  + q_\text R^2 q_\text L + q_\text R^3 }{4}
\end{align}
E.g. a shock joining 1 and $-5$ is moving at speed $-26$.

\begin{figure}[h]
 \centering
 \includegraphics[width=\textwidth]{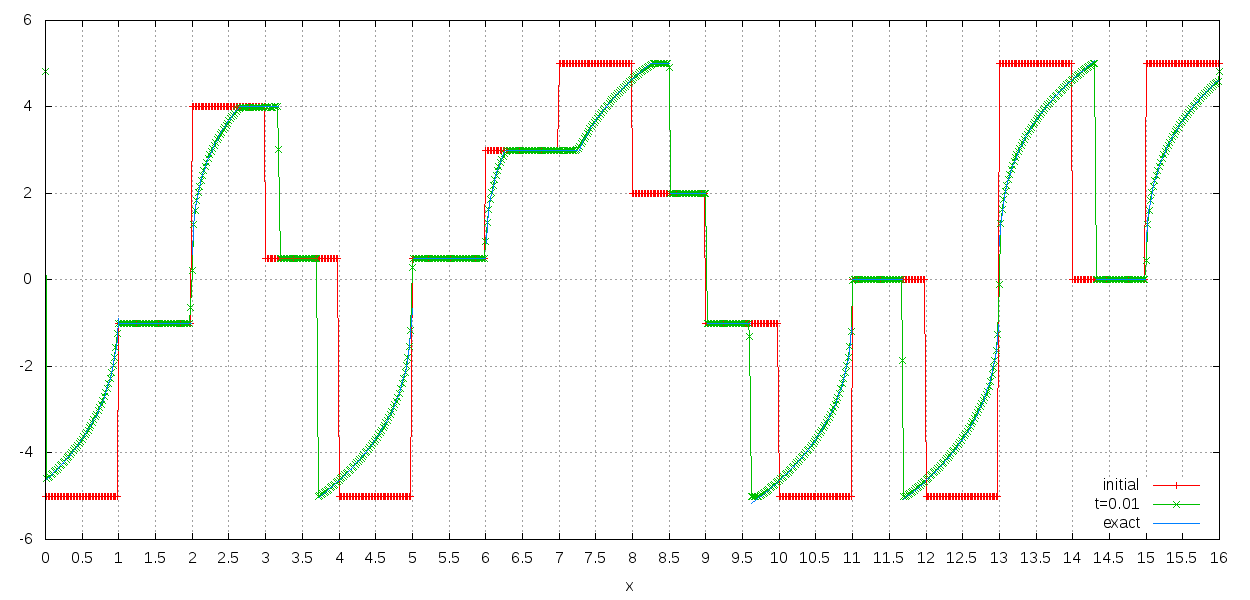}
 \caption{Equation \eqref{eq:quarticequation} evolved with the Active Flux scheme using the modified fixpoint iteration \eqref{eq:selectingchar}. All the shock speeds and rarefactions are correctly resolved. The solid line shows the exact solution. The initial data are piecewise constant. Point values of the solution are shown at time $t = 0.01$. The boundaries are periodic. $\Delta x = 2/100$ and the CFL number is 0.45. Limiting procedure according to Theorem \ref{thm:powerlawlimiting} / Equation \eqref{eq:limitingonepowerlaw} is used.}
 \label{fig:cubicallriemann} 
\end{figure}

The self-similar rarefaction is given by $\sqrt[3]{\frac{x}{t}}$. Fig. \ref{fig:cubicallriemann} shows again a selection of shocks and rarefaction together with the numerical and the exact solution. Usage of the modified fixpoint iteration \eqref{eq:selectingchar} allows to resolve all the shocks and rarefactions correctly.

\newpage
\subsection{$p$-system} \label{ssec:psystem}

The $p$-system
\begin{align}
 \del_t \rho + \del_x v &= 0\\
 \del_t v + \del_x p(\rho) &= 0 & p(\rho) &= \rho^\gamma
\end{align}
in case of smooth solutions can be rewritten in the form \eqref{eq:nonlinsysdiagonal} as
\begin{align}
 \del_t \left(\frac{2 \rho c}{\gamma + 1} \pm v \right ) \pm c \del_x \left(\frac{2 \rho c}{\gamma + 1} \pm v \right ) &= 0 
\end{align}
with $c = \sqrt{p'(\rho)}$. In the following, $\gamma = 1.4$ is used.

As the eigenvalues cannot switch sign, and never vanish, Active Flux solving the $p$-system seems less prone to artefacts.

\subsubsection{Convergence study}

\begin{figure}[h]
 \centering
 \includegraphics[width=0.48\textwidth]{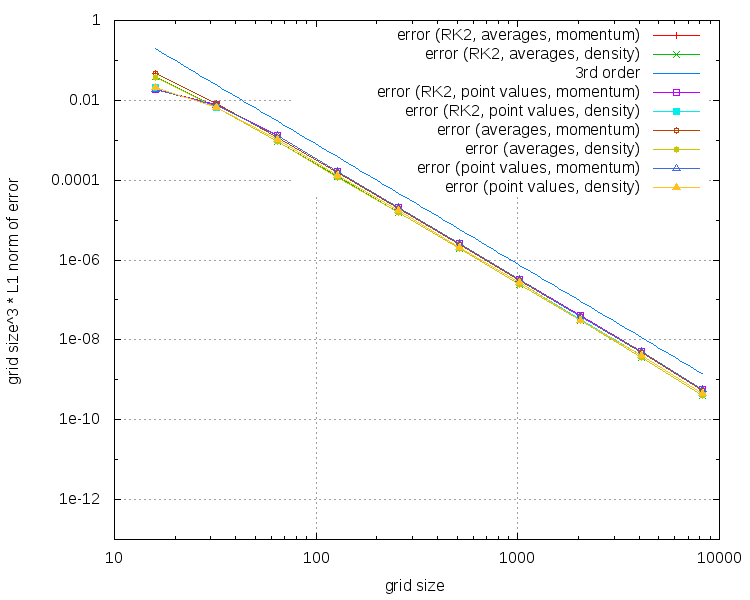}\hfill\includegraphics[width=0.48\textwidth]{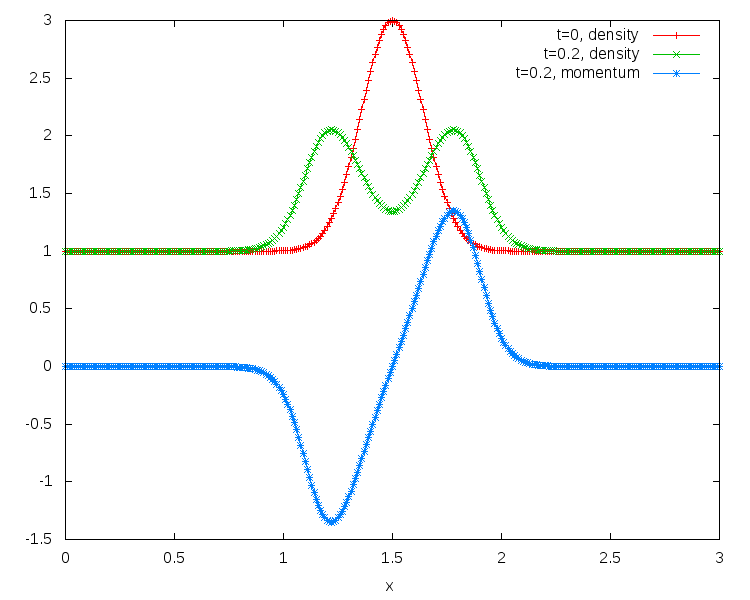}
 \caption{The $p$-system evolved with the Active Flux scheme. \textit{Left}: Third order convergence of the numerical solution on both point values and averages, for momentum $v$ and density $\rho$, using both the algorithm from section \ref{ssec:systemsaveragespeed} and from section \ref{ssec:rungekutta} (the latter marked \texttt{RK2}). The lines showing the results for different schemes and quantities virtually lie on top of each other indicating comparable error. \textit{Right}: Setup and numerical solution for $\Delta x = 1/100$ showing point values. No limiting used.  
 }
 \label{fig:psystemconvergence}
\end{figure}

Figure \ref{fig:psystemconvergence} (left) demonstrates the correct order of convergence for both algorithms of section \ref{ssec:systemsaveragespeed} and \ref{ssec:rungekutta}. The initial setup (shown in Figure \ref{fig:psystemconvergence} (right)) is a Gaussian in the density and $v\equiv 0$. During the evolution two waves are forming which self-steepen. At the final time $t=0.2$ they have not formed shocks, though. The reference solution is obtained by evolving the problem on a highly resolved grid of $16384=2^{14}$ cells with the algorithm of section \ref{ssec:rungekutta}.

\subsubsection{Riemann problem}

Consider the following Riemann problem
\begin{align}
 v = \begin{cases} \rho = 2, v = 1 & 0.3< x <0.7 \\ \rho = 0.1, v = -0.5  & \text{else} \end{cases}
\end{align}

\begin{figure}[h]
 \centering
 \includegraphics[width=0.48\textwidth]{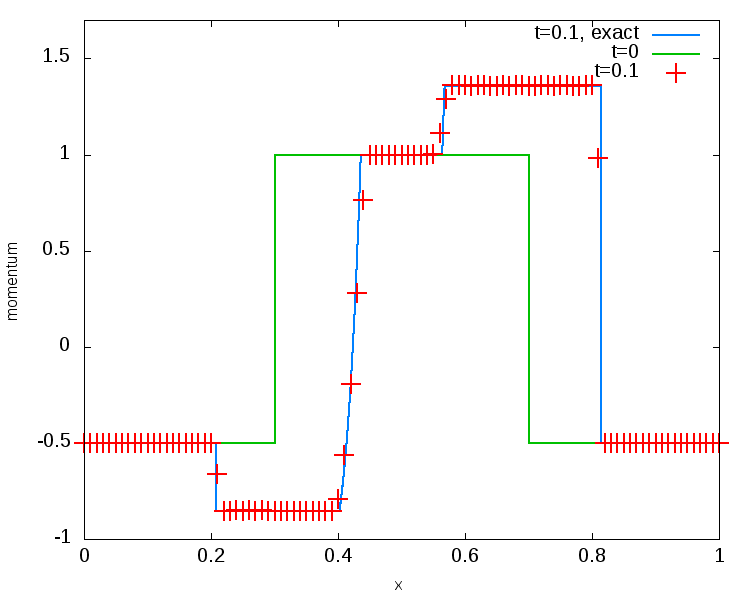} \hfill \includegraphics[width=0.48\textwidth]{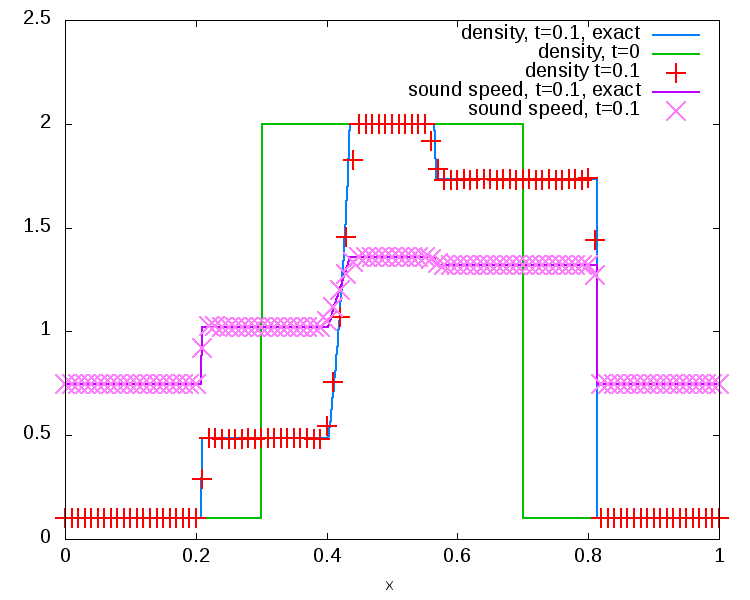}
 \caption{The $p$-system evolved with the Active Flux scheme using the iteration from section \ref{ssec:rungekutta} (no difference was apparent upon using the iteration from \ref{ssec:systemsaveragespeed}). Here, $\Delta x = 1/100$ and CFL = 0.45; the limiting function \eqref{eq:limitingsymmetrizedpowerlaw} is applied to the conserved quantities. Periodic boundaries are used and the curves show the point values at time $t=0.1$. \textit{Left}: Momentum $v$. \textit{Right}: Density $\rho$ and speed of sound $c$.} 
 \label{fig:psystemrp} %
\end{figure}

\subsection{Isentropic Euler equations} \label{ssec:euler}

Consider the isentropic Euler equations
\begin{align}
 \del_t \rho + \del_x (\rho v) &= 0\\
 \del_t (\rho v) + \del_x (\rho v^2 + p(\rho)) &= 0 & p(\rho) &= K \rho^\gamma
\end{align}
On smooth solutions this system is equivalent to
\begin{align}
\del_t \left( v \pm \frac{2c}{\gamma-1}    \right ) + (v \pm c) \del_x \left(v \pm \frac{2c}{\gamma-1}\right) &= 0 
\end{align}
with $c = \sqrt{\gamma p(\rho) / \rho}$. In the following, $K=1$, $\gamma=1.4$ are used.

\subsubsection{Convergence study}

\begin{figure}[h]
 \centering
 \includegraphics[width=0.48\textwidth]{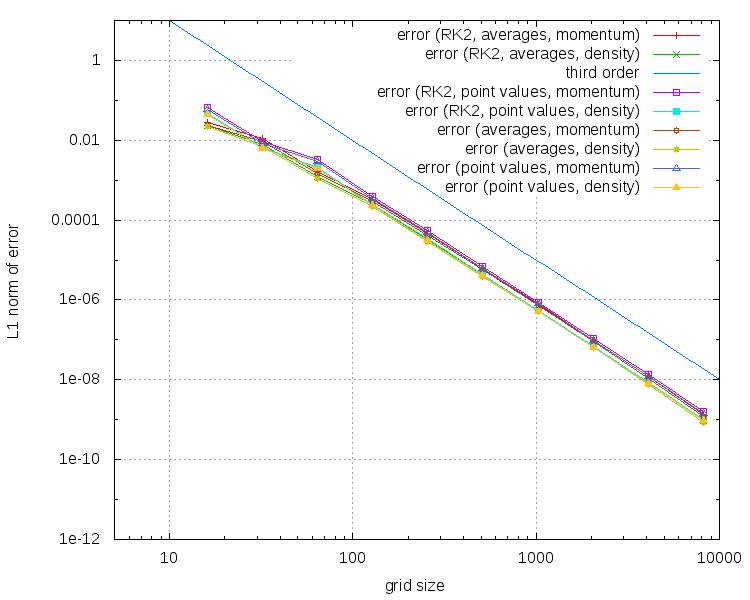}\hfill\includegraphics[width=0.48\textwidth]{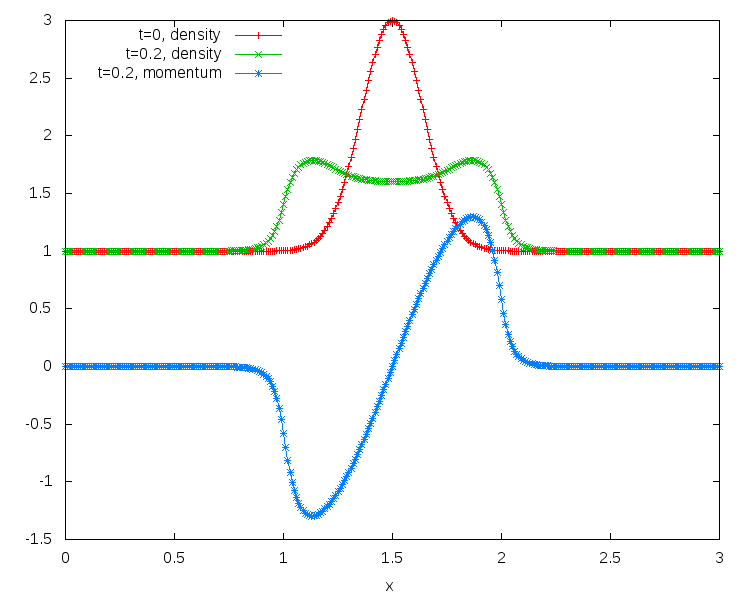}
 \caption{Isentropic Euler equations solved with the Active Flux scheme. \textit{Left}: Third order convergence of the numerical solution on both point values and averages, for momentum $\rho v$ and density $\rho$, using both the algorithm from section \ref{ssec:systemsaveragespeed} and from section \ref{ssec:rungekutta} (the latter marked \texttt{RK2}). The lines showing the results for different schemes and quantities virtually lie on top of each other indicating comparable error. \textit{Right}: Setup and numerical solution for $\Delta x = 1/100$ showing point values. No limiting used.}
 \label{fig:eulerconvergence}  
\end{figure}

Figure \ref{fig:eulerconvergence} shows the setup of a Gaussian initial density and no velocity. Its evolution (before a shock forms) is used to study the convergence of the method. The reference solution is the setup solved with the iteration from section \ref{ssec:rungekutta} on a grid of $16384$ cells. Third order is confirmed experimentally.

 \newpage
 
\subsubsection{Riemann problems}

\begin{figure}[h]
 \centering
 \includegraphics[width=0.48\textwidth]{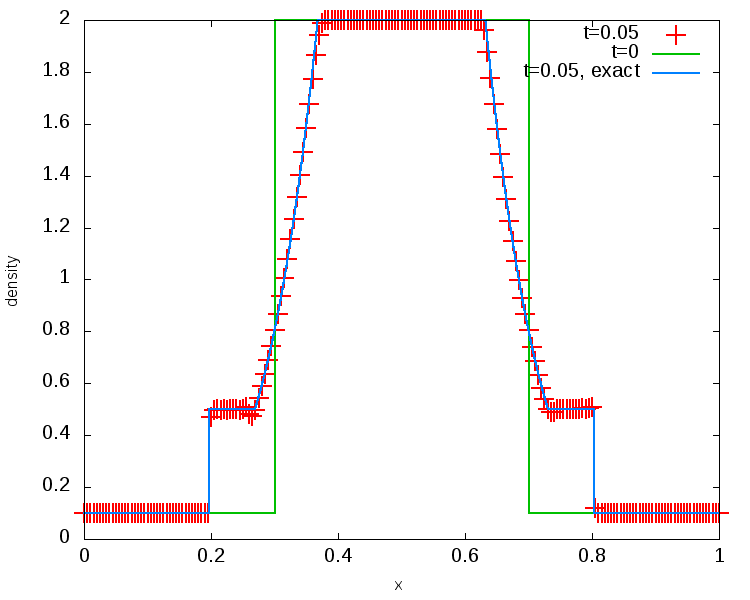} \hfill \includegraphics[width=0.48\textwidth]{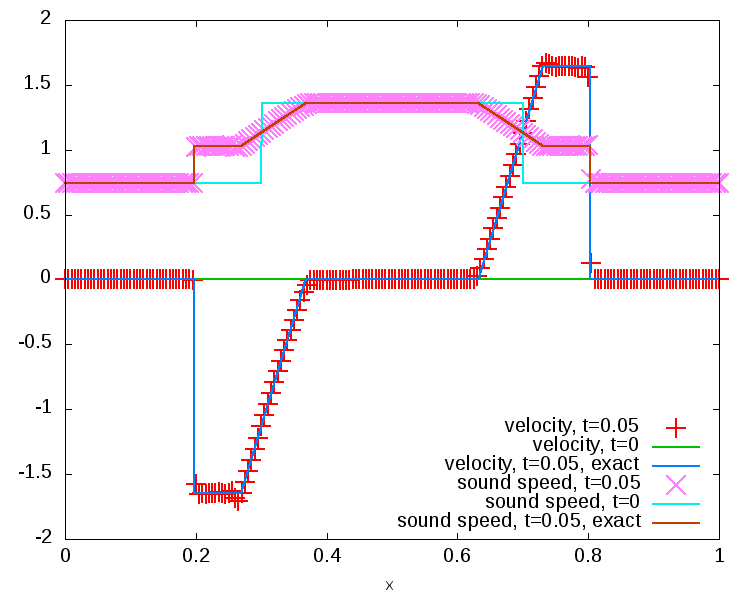}
 \caption{Isentropic Euler equations solved with the Active Flux scheme using iteration from section \ref{ssec:rungekutta} and the fix \eqref{eq:rungekuttastep1modif}--\eqref{eq:rungekuttamodifselect}. Here, $\Delta x = 1/200$ and CFL = 0.45; the limiting function \eqref{eq:limitingonepowerlaw} is applied to the conserved quantities. Periodic boundaries are used and the curves show the point values at time $t=0.05$. \textit{Left}: Density $\rho$. \textit{Right}: Velocity $v$, sound speed $c$.} 
 \label{fig:psystemrp1} %
\end{figure}

\begin{figure}[h]
 \centering
 \includegraphics[width=0.48\textwidth]{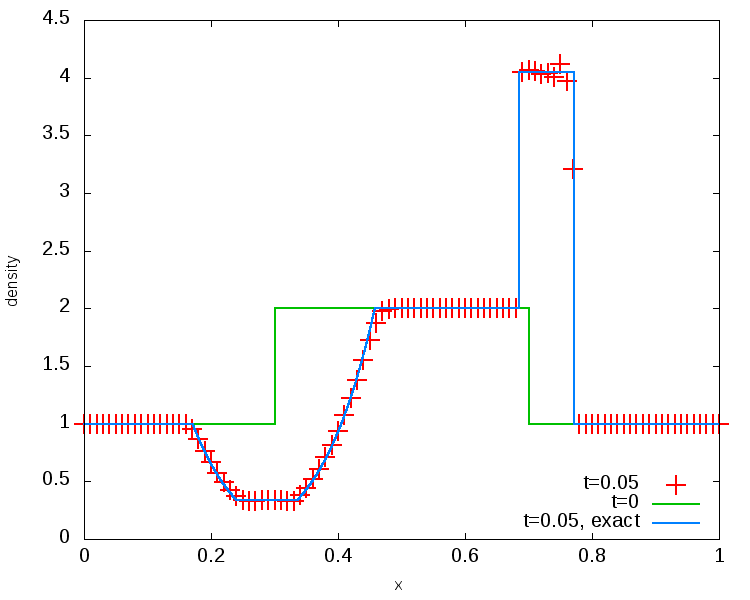} \hfill \includegraphics[width=0.48\textwidth]{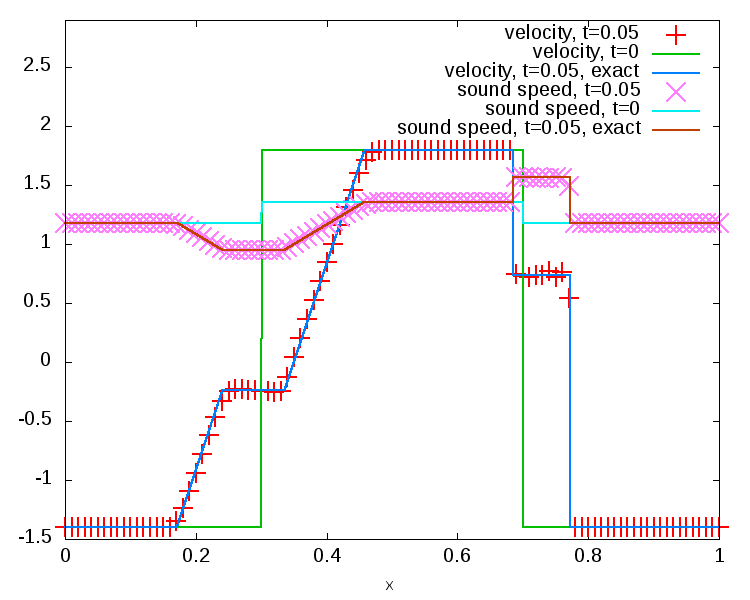}
 \caption{Isentropic Euler equations solved with the Active Flux scheme. Apart from the values for the states of the Riemann problem ($\rho \in \{1, 2 \}, v \in \{-1.4, 1.8 \}$), solution method and plot as in Figure \ref{fig:psystemrp1}.} 
 \label{fig:psystemrp2} %
\end{figure}

The scheme is now applied to Riemann problems. Here it becomes important to use the modification \eqref{eq:rungekuttastep1modif}--\eqref{eq:rungekuttamodifselect}. Figure \ref{fig:psystemrp1} shows a setup with a transonic rarefaction, accurately resolved, and Figure \ref{fig:psystemrp2} shows a strong shock and a double rarefaction. The two tests cover the cases of the two eigenvalues of the system having different and same sign. The double shock case shows little artefacts which vanish upon refining the grid.

\clearpage

\subsection{Full Euler equations} \label{ssec:fulleuler}

The full Euler equations
\begin{align}
 \del_t \rho + \del_x (\rho v) &= 0\\
 \del_t (\rho v) + \del_x (\rho v^2 + p) &= 0\\
 \del_t e + \del_x (v (e+p)) &= 0
\end{align}
with the ideal equation of state
\begin{align}
 e &= \frac{p}{\gamma-1} + \frac12 \rho v^2 & \gamma &> 1
\end{align}
form a hyperbolic system of conservation laws, but do not admit characteristic variables. Thus, the most general solution algorithm \eqref{eq:approxevoopgeneral} is required. In the following, $\gamma = 1.4$ is used.

\subsubsection{Convergence study}

\begin{figure}[h]
 \centering
 \includegraphics[width=0.48\textwidth]{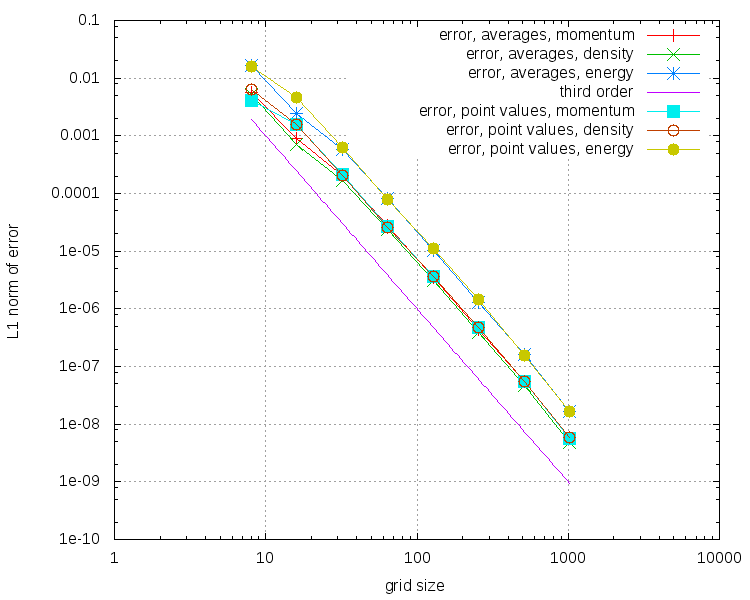}\hfill\includegraphics[width=0.48\textwidth]{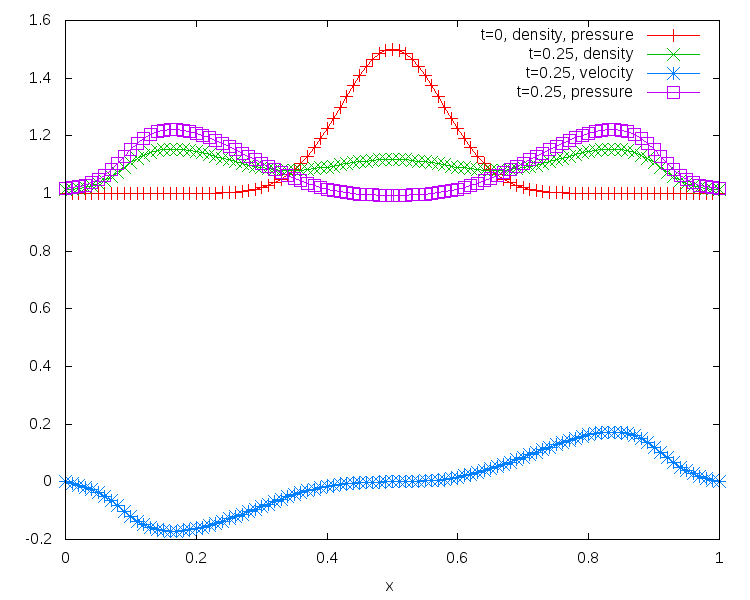}
 \caption{Full Euler equations solved with the Active Flux scheme using the approximate evolution operator \eqref{eq:approxevoopgeneral}. \textit{Left}: Third order convergence of the numerical solution on both point values and averages, for momentum $\rho v$, density $\rho$ and energy $e$. The lines virtually lie on top of each other indicating comparable error. A CFL number of 0.7 is used. \textit{Right}: Setup and numerical solution for $\Delta x = 1/100$ showing point values. No limiting used.}
 \label{fig:fulleulerconvergence}  
\end{figure}

To demonstrate the convergence of the method, the test from \cite{kerkmann18} is run:
\begin{align}
 \rho_0(x) &= p_0(x) = 1 + \frac12 e^{-80 (x-0.5)^2} & v_0(x) &= 0
\end{align}
The numerical results are compared at $t=0.25$ to a reference solution obtained on a grid of 2048 points. In Figure \ref{fig:fulleulerconvergence}, one observes third order convergence, in agreement with the theoretical expectation.

\subsubsection{Riemann problem}

\begin{figure}
 \centering
 \includegraphics[width=0.45\textwidth]{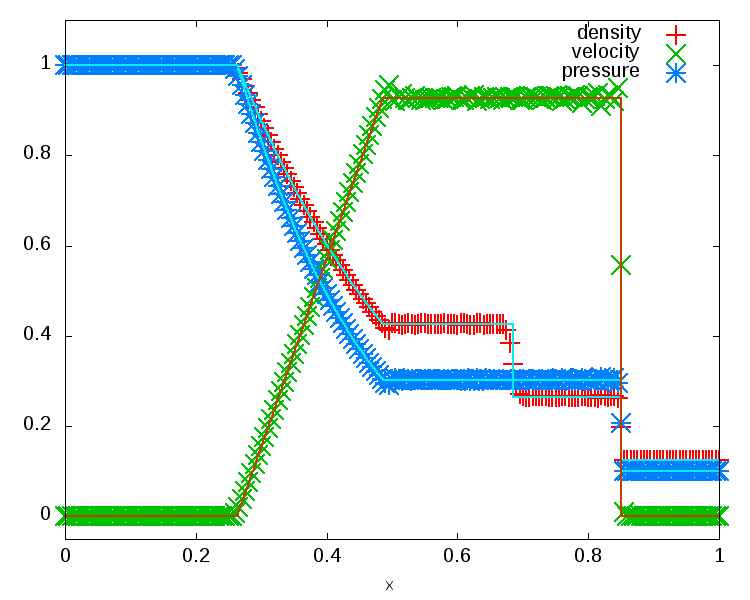} \hfill \includegraphics[width=0.45\textwidth]{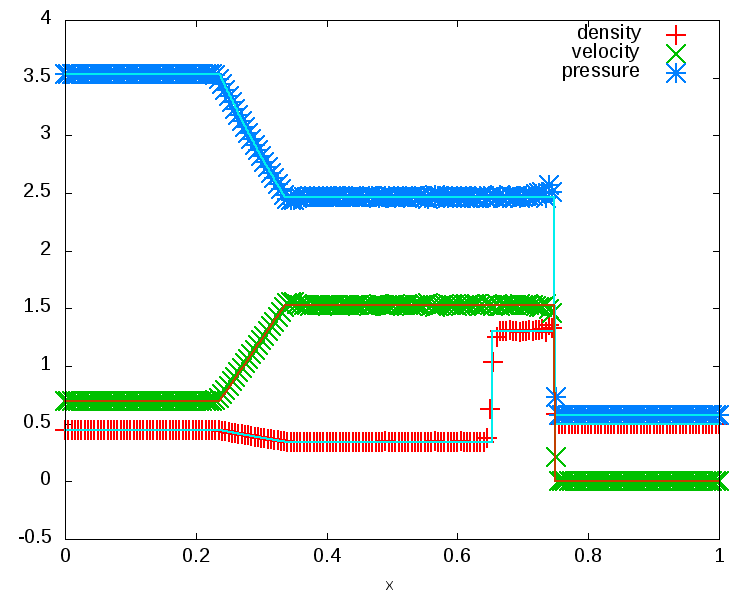}
 \caption{Riemann problem setups for the full Euler equations. The approximate solution operator \eqref{eq:approxevoopgeneral} with limiting is used on a grid with $\Delta x = 1/200$ and with a CFL number of 0.7. Point values are shown at $t=0.1$. \emph{Left}: Sod's test problem (\cite{sod78}). \emph{Right}: Lax's test problem (\cite{lax54}). Solid lines show the exact solution.}
 \label{fig:fulleulerriemann}
\end{figure}

To assess the performance of the numerical method on discontinuous problems, two Riemann problems are computed: the Sod shock tube (\cite{sod78}, Figure \ref{fig:fulleulerriemann}, \emph{left}) and the Lax shock tube (\cite{lax54}, Figure \ref{fig:fulleulerriemann}, \emph{right}). One observes the good perfromance of Active Flux even on these discontinuous setups.

\subsubsection{Interaction between a shock and a sound wave}

\begin{figure}
 \centering
 \includegraphics[width=0.45\textwidth]{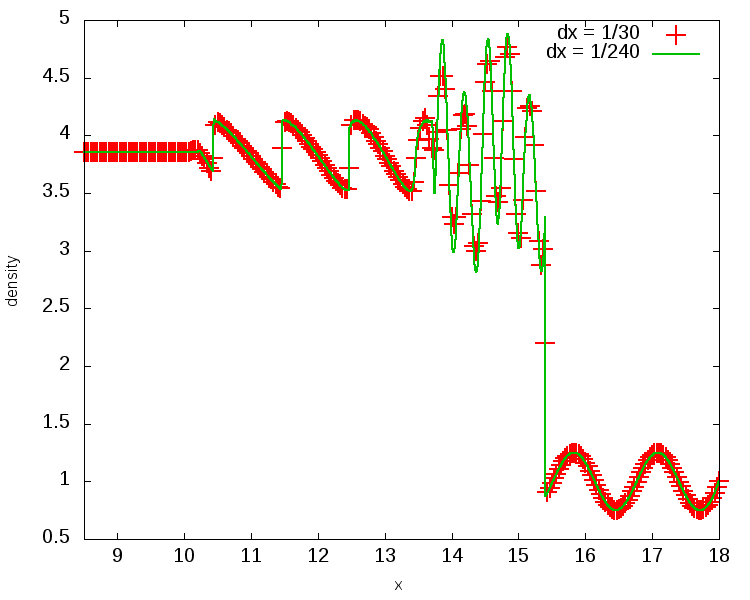} \hfill \includegraphics[width=0.45\textwidth]{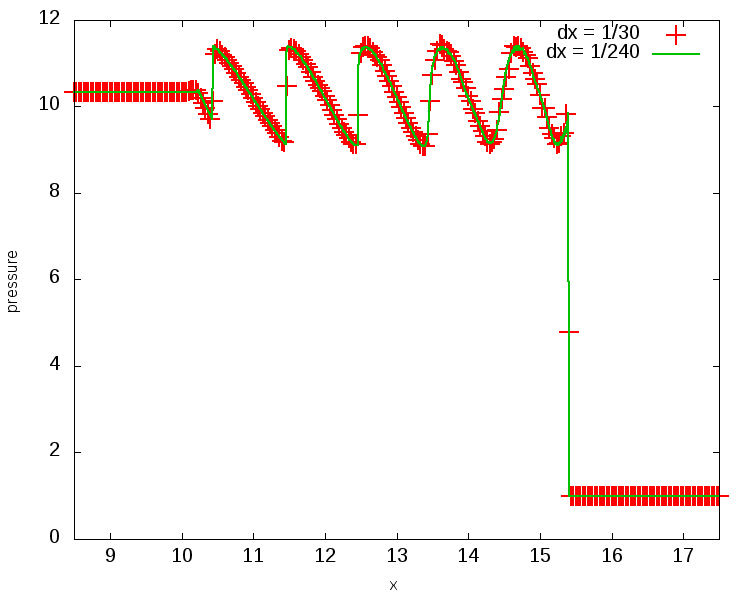}
 \caption{The Active Flux method is used to solve the Shu-Osher test (\cite{shu89}). The approximate solution operator \eqref{eq:approxevoopgeneral} with limiting is used with a CFL number of 0.7 and on grids with $\Delta x = 1/30$ (crosses) and $1/240$ (solid line). Point values are shown at $t=0.18$. \emph{Left}: Density. \emph{Right}: Pressure.}
 \label{fig:fulleulershuosher}
\end{figure}

Finally, to demonstrate the performance of the algorithm on a more challenging setup, the Shu-Osher test (\cite{shu89}) is shown in Figure \ref{fig:fulleulershuosher}. One observes that due to its high order, Active Flux is able to capture the details of the interaction on poorly resolved grids without difficulty (compare e.g. to \cite{cockburn89}).

\subsection{Multi-dimensional scalar equations} \label{ssec:scalarnumericsmultid}

As a last test case, consider the multi-dimensional Burgers' equation
\begin{align}
 \del_t q + \del_x \left( \frac{q^2}{2} \right ) + \del_y \left( \frac{q^2}{2}   \right ) = 0 \label{eq:multidburgerseq}
\end{align}
and a 4-quadrant Riemann problem setup as follows:
\begin{align}
 q_0 = \begin{cases} -1 & \text{NE} \\ -0.2 & \text{NW} \\ 0.5 & \text{SW} \\ 0.8 & \text{SE} \end{cases}
\end{align}

\begin{figure}[h]
 \centering
  \includegraphics[width=0.48\textwidth]{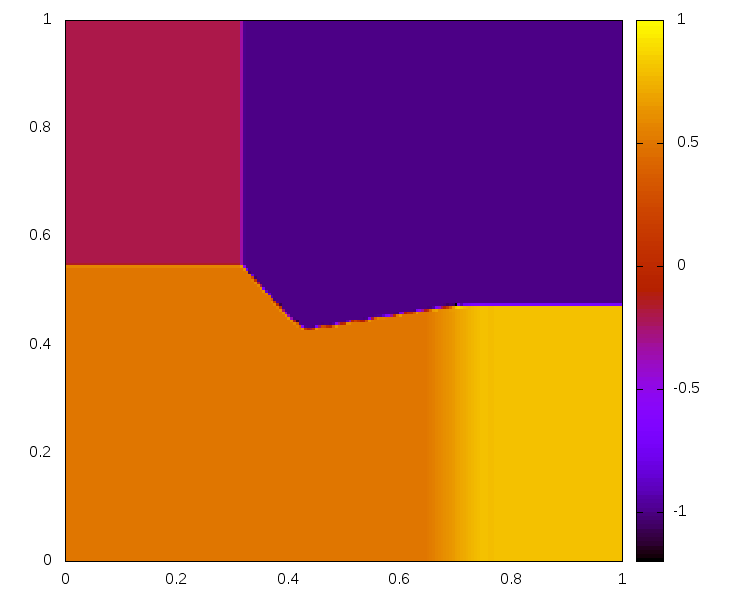}\hfill \includegraphics[width=0.48\textwidth]{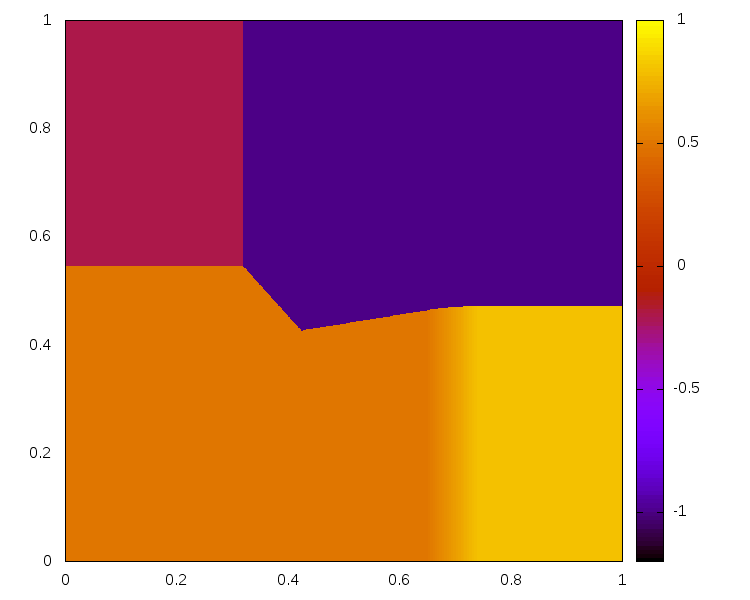}
  \caption{The multi-dimensional Burgers' equation \eqref{eq:multidburgerseq} solved with the Active Flux scheme using the fixpoint iteration \eqref{eq:selectingchar}. Here, $\Delta x = \Delta y = 1/200$ and CFL = 0.9. The solution has been computed on a grid of double size to avoid the influence of the boundaries. \textit{Left}: Cell averages of the solution at time $t=0.3$. \textit{Right}: Exact solution following \cite{guermond11}.}
  \label{fig:multidburgers}
\end{figure}

Figure \ref{fig:multidburgers} shows the solution at $t=0.3$ using Active Flux (with fixpoint iteration \eqref{eq:selectingchar}) along with the exact solution taken from \cite{guermond11} (p. 4258). No limiting is used.

\clearpage

\section{Conclusion and outlook}

The Active Flux scheme is a finite volume scheme with additional pointwise degrees of freedom located at the cell boundary. It involves a continuous reconstruction and thus does not make use of Riemann solvers. Instead, an evolution operator for the pointwise degrees of freedom is required. Once their evolution is obtained, the update of the cell average follows the usual finite volume/\-Godunov scheme idea; the intercell flux is obtained by evaluating the flux function on the point values along the boundary (and using quadrature). The Active Flux scheme has initially (\cite{vanleer77,eymann13}) been employing exact evolution operators, because the problems onto which the scheme was applied admitted closed form solution operators. 

In particular it has been shown in \cite{barsukow18activeflux} that for linear acoustics the Active Flux method is low Mach number compliant and stationarity preserving without the need for any fix. This makes Active Flux an interesting candidate for a class of methods, which are structure preserving by construction. Inspired by the finding for linear acoustics, this paper serves as a stepping stone towards deriving structure preserving Active Flux methods for multi-dimensional nonlinear systems.

Upon an extension of the Active Flux method to nonlinear problems usage of exact evolution operators cannot be maintained. Approximate evolution operators suggested so far in the literature either were reducing the order of the scheme or involved complicated expressions such as the Lax-Wendroff expansion and subsequent solution of Riemann problems (ADER). This paper shows how approximate evolution operators can be found which are not costly and allow to maintain third order of accuracy. The cases considered here are nonlinear scalar conservation laws in one and two spatial dimensions as well as nonlinear hyperbolic systems of conservation laws in one spatial dimension. 

As Active Flux is very different from standard finite volume or Galerkin methods, many aspects need to be reconsidered, and many questions well-studied for other methods were still open. Continuous reconstruction might raise doubts about the applicability of Active Flux to problems involving shock formation. 
It is found that in certain setups too simple an evolution operator fails to correctly recognize the self-steepening. This then leads to artefacts that resemble entropy glitches encountered with certain finite volume schemes. As a cure, in this paper a simple strategy is presented which allows to take into account the fact that characteristics may cross. This ``entropy fix'' is found to lead to accurate numerical evolutions without artefacts. By means of numerical examples it is shown that Active Flux, for example, is able to accurately solve Riemann problems for one-dimensional systems of nonlinear conservation laws, such as the Euler equations.

As a numerical method of higher order is prone to overshoots around discontinuities, a limiting procedure needs to be in place. Here, a simple limiting is suggested which modifies the reconstruction whenever an avoidable overshoot/undershoot is recognized. Whereas the Active Flux limiters available in the literature either require joining several polynomials inside the cell or re-introduce discontinuities, the suggested monotone reconstruction is simple to compute and retains a continuous reconstruction.

The correct approximation of the entropy solution and limiting in one spatial dimension may not outperform currently available methods of third and higher order. However, all these are necessary ingredients for an extension to multiple spatial dimensions that so far were open, or at least insufficiently studied for the Active Flux method.

Future work shall be devoted to multi-dimensional hyperbolic systems.
The approximate evolution operators presented here shall be extended to multiple spatial dimensions and thus combined with the favorable properties of Active Flux in multiple dimensions. This hopefully will pave the way towards a powerful structure-preserving method for multi-dimensional systems of conservation laws.

\section*{Acknowledgement}

The author was supported by the German Academic Exchange Service (DAAD) with funds from the German Federal Ministry of Education and Research (BMBF) and the European Union (FP7-PEOPLE-2013-COFUND -- grant agreement no. 605728), as well as by the Deutsche Forschungsgemeinschaft (DFG) through project 429491391.

\end{document}